\documentclass[12pt]{article}
\usepackage{graphics,graphicx}
\usepackage{amssymb,amsmath,amsopn,amsfonts}
\usepackage[dvips]{color}
\usepackage{colordvi,multicol}
\usepackage{epsf}
\newfam\msbfam
\font\tenmsb=msbm10 \textfont\msbfam=\tenmsb \font\sevenmsb=msbm7
\scriptfont\msbfam=\sevenmsb \font\fivemsb=msbm5
\scriptscriptfont\msbfam=\fivemsb

\def\th#1{\vspace{1mm}\noindent{\bf #1}\quad}

\def\proof{\vspace{1mm}\noindent{\it Proof}\quad}

\numberwithin{equation}{section}

\def\bc{\begin{center}}
\def\ec{\end{center}}
\def\no{\noindent}
\def\hang{\hangindent\parindent}
\def\textindent#1{\indent\llap{\qquad #1\ \ \enspace}\ignorespaces}
\def\ref{\par\hang\textindent}

\begin{document}

\title{ {\bf Dirichlet form associated with the  $\Phi_3^4$ model
\thanks{Research supported in part  by  NSFC ( No.11401019, No.11671035) and DFG through  CRC
701
}\\} }
\author{  {\bf Rongchan Zhu}$^{\mbox{a,c}}$, {\bf Xiangchan Zhu}$^{\mbox{b,c},}$\thanks{Corresponding author}
\date{}
\thanks{E-mail address:
zhurongchan@126.com(R. C. Zhu), zhuxiangchan@126.com(X. C. Zhu)}\\\\
$^{\mbox{a}}$Department of Mathematics, Beijing Institute of Technology, Beijing 100081,  China\\
$^{\mbox{b}}$School of Science, Beijing Jiaotong University, Beijing 100044, China\\
$^{\mbox{c}}$Department of Mathematics, University of Bielefeld, D-33615 Bielefeld, Germany}

\maketitle

\noindent {\bf Abstract}
We construct the Dirichlet form associated with the dynamical $\Phi^4_3$ model obtained in [Hai14, CC13] and [MW16].
This Dirichlet form on cylinder functions is identified as a classical gradient bilinear form.
As a consequence,  this classical gradient bilinear form is closable and then by a well-known result its closure is also a quasi-regular Dirichlet form, which means that there exists another (Markov) diffusion process, which also admits the $\Phi^4_3$ field measure as an invariant (even symmetrizing) measure.

\vspace{1mm}
\no{\footnotesize{\bf 2000 Mathematics Subject Classification AMS}:\hspace{2mm} 60H15, 82C28}
 \vspace{2mm}

\no{\footnotesize{\bf Keywords}:\hspace{2mm}  $\Phi_3^4$ model, Dirichlet form, regularity structures, paracontrolled distributions, space-time white noise, renormalisation}

\section{Introduction}

Recall that the usual continuum Euclidean $\Phi^4_d$-quantum field theory is heuristically described by the following
probability measure:
$$\mu(dx)=N^{-1}\Pi_{\xi\in\mathbb{T}^d}dx(\xi)\exp\bigg(-\int_{\mathbb{T}^d}(|\nabla x(\xi)|^2+mx^2(\xi)+\frac{\lambda}{2}x^4(\xi))d\xi\bigg),\eqno(1.1)$$
where $N$
is the normalization constant, $m$ is a real constant, $\lambda\geq0$ is the coupling constant and $x$ is the real-valued field and $\mathbb{T}^d$ is the $d$-dimensional torus. There have been many approaches to the problem of
giving a meaning to the above heuristic measure for $d=2$ and $d=3$ (see  [GRS75] [GJ87] and references
therein). The construction of this $\Phi^4_3$ field measure $\mu$ has been achieved in [Fel74] for $\lambda$ small enough, which was one of the major achievements
of the programme of constructive quantum field theory.
In [PW81]  Parisi and Wu proposed a program for Euclidean quantum field
theory of getting Gibbs states of classical
statistical mechanics as limiting distributions of stochastic processes,  especially as solutions to non-linear stochastic differential
equations. Then one can use the stochastic differential equations to study the properties of the Gibbs states. This
procedure is called stochastic field quantization (see [JLM85]). The $\Phi_d^4$ model is the simplest non-trivial Euclidean quantum field (see [GJ87] and the reference therein). The issue of the stochastic quantization of the $\Phi^4_d$ model is to solve the following equation:
$$\aligned d \Phi= &( \Delta \Phi-\lambda\Phi^3-m\Phi)dt +dW(t) \quad\Phi(0)=\Phi_0.\endaligned\eqno(1.2)$$
where $W$ is a cylindrical Wiener process on $L^2(\mathbb{T}^d)$. In the following we take $\lambda$ small enough (weak coupling) as in [BFS83] and in the following when we analyze (1.2) we  omit $\lambda$  for simplicity if there is no confusion. The solution $\Phi$ is also called  dynamical $\Phi^4_d$ model. The main difficulty in this case is that $W$
and hence the solutions $\Phi$ are so singular that the non-linear term is not well-defined in the classical sense.

In
two spatial dimensions, the dynamical $\Phi_2^4$ model was first treated in [AR91] by using the Dirichlet form approach: The authors considered the following bilinear form on  $L^2(E;\mu)$ with $E$ being a separable Banach space and $\mu(E)=1$:
$$\mathcal{E}(u,v):=\frac{1}{2}\int \langle Du, Dv\rangle_{L^2}d\mu,$$
where $Du$ means $L^2$-derivative, which is defined in Section 4. By the corresponding integration by parts formula for $\mu$ they obtained that the bilinear form is closable and its closure $(\mathcal{E},D(\mathcal{E}))$ is a quasi-regular Dirichlet form. Then according to a general result in [MR92] (see Theorem D.4), we know that there exists a (Markov) diffusion process $M=(\Omega,\mathcal{F},X(t),$ $(P^x)_{x\in E})$ on $E$ \emph{properly associated with} $(\mathcal{E},D(\mathcal{E}))$. The sample paths of the associated  process  satisfy (1.2) in the (probabilistically) weak sense for  quasi-surely every $\Phi_0$.

Later in [DD03] and [MW15], the authors split $\Phi$ as $\Phi=\Phi_1+v$, where
$$d\Phi_1=\Delta \Phi_1dt+dW,$$
$$\partial_t v=\Delta v-(v^3+3v^2\Phi_1+3v:\Phi_1^2:+:\Phi_1^3:)-m(\Phi_1+v),\eqno(1.3)$$
where $:\Phi_1^2:, :\Phi_1^3:$ are defined as Wick products. Then the nonlinear terms are well defined in the classical sense and they obtained a (probabilistically) strong solution to (1.3).

In three
spatial dimensions both techniques  break down. For the Dirichlet form approach we cannot directly obtain that
 the  bilinear form:
$$\mathcal{E}(u,v):=\frac{1}{2}\int_E \langle Du, Dv\rangle_{L^2}d\mu, $$
is closable since the measure $\mu$ is more singular and may be not quasi-invariant along smooth direction (see [ALZ06]). Nobody has constructed the Dirichlet form associated with $\Phi^4_3$ model successfully and the closablity of the corresponding bilinear form has been a long-standing open problem for more than 25 years ([AR91]).  For the second approach (1.3) is also not well defined in the classical sense since the noise is more rough.
It was a long-standing open problem to give a meaning to the equation (1.2) in the three dimensional
case. A breakthrough result was achieved recently by Martin Hairer in [Hai14], where he introduced
a theory of regularity structures and gave a meaning to  equation (1.2) successfully. Also by using the paracontrolled distributions proposed by Gubinelli, Imkeller and Perkowski
in [GIP15] existence and uniqueness of local solutions to (1.2) have been obtained in [CC13].
Recently, these two approaches have been successful in giving a meaning to a lot of ill-posed
stochastic PDEs like the Kardar-Parisi-Zhang (KPZ) equation ([KPZ86], [BG97], [Hai13]), the stochastic 3D-Navier-Stokes equation driven by space-time white
noise ([ZZ14], [ZZ15a]), the dynamical sine-Gordon equation ([HS16]) and so on (see [HP14]
for more other interesting examples). These two approaches  are inspired by the theory
of rough paths [Lyo98]. In [Kup16] the author also uses renormalization group
techniques to make sense of the dynamical  $\Phi^4_3$ model. Recently in [MW16] the authors obtained  global well-posedness of the solution to (1.2) in the three dimensional case based on the paracontrolled distribution method.

The aim of  this paper is to construct the Dirichlet form associated to the $\Phi_3^4$ model. Dirichlet form techniques have developed into a powerful
method to combine analytic and functional analysis, as well as potential
theoretic and probabilistic methods to study the properties of
stochastic processes. In [RZZ15, RZZ16]  M. R\"{o}ckner and the authors of this paper combine the Dirichlet form approach and the SPDE approach to  obtain new properties in the two dimensional case (such as restricted Markov uniqueness and the characterization of the $\Phi^4_2$ field).  We hope this paper is a start to study the dynamical $\Phi^4_3$ model combining  Dirichlet form techniques and the theory of regularity structures as well as the paracontrolled distributions approach.

  Different from [AR91], our idea is to construct the Dirichlet form from the global solution $\Phi(t)$ obtained in [MW16]. It has been proved in [HM15] that $\Phi(t)$ satisfies Markov property. Moreover, it is easy to obtain that $\Phi(t)$ satisfies the Feller property (see Lemma 4.1), which implies that $\Phi(t)$  satisfies the  strong Markov property. Then we  prove $\Phi(t)$ is reversible with respect to $\mu$ by the lattice approximations obtained in [ZZ15] (see Lemma 4.2). Hence we obtain our first main result of this paper:
\vskip.10in
\th{Theorem 1.1} There exists a quasi-regular Dirichlet form $(\mathcal{E}, D(\mathcal{E}))$ associated with $\Phi(t)$. Moreover, $\Phi$ is \emph{properly associated with} $(\mathcal{E}, D(\mathcal{E}))$ in the sense that the semigroup for $\Phi$ is a quasi-continuous version of the semigroup associated with $(\mathcal{E}, D(\mathcal{E}))$. Furthermore, $\mathcal{F}C_b^\infty\subset D(\mathcal{E})$ and $\langle l,\cdot\rangle\in D(\mathcal{E})$ for any $l\in E^*$.
\vskip.10in
For definitions of quasi-regular Dirichlet form we refer to Appendix D. Here  $\mathcal{F}C_b^\infty$ denotes all the smooth with all derivatives bounded cylinder functions on the state space $E$, $E^*$ is the dual space of $E$ and $\langle \cdot,\cdot\rangle$ is the dualization between $E$ and $E^*$. For the explicit definition we refer to Section 4. Moreover, we can identify the Dirichlet form on the cylinder functions as a gradient Dirichlet form:

\vskip.10in
\th{Theorem 1.2} For $f,g\in \mathcal{F}C_b^\infty$, $\mathcal{E}(f,g)=\frac{1}{2}\int \langle Df,Dg\rangle d\mu$ with $\langle\cdot,\cdot\rangle$ being the inner product of $L^2(\mathbb{T}^3)$ and $Df$ is $L^2$-derivative defined in Section 4.
\vskip.10in
 As a byproduct of Theorem 1.2 we can also deduce that $\Phi$ is an energy solution in the stationary case (see Remark 5.2). Energy solution is a notion of weak solutions for KPZ equation to describe the large scale fluctuations of a wide class of
weakly asymmetric particle systems (see [GJ13, GJ13a, GP15]). For the dynamical $\Phi^4_3$ case we can also introduce the notion of energy solution.

As a consequence of Theorem 1.2, we obtain that the bilinear form is closable, which we cannot directly obtain as we mentioned before:
\vskip.10in
\th{Theorem 1.3} The bilinear form $\bar{\mathcal{E}}(f,g)=\frac{1}{2}\int \langle Df, Dg\rangle d\mu, f,g\in\mathcal{F}C_b^\infty$, is closable and its closure $(\bar{\mathcal{E}},D(\bar{\mathcal{E}}))$ is a quasi-regular Dirichlet form. Then  there exists a (Markov) diffusion process  \emph{properly associated with} $(\bar{\mathcal{E}},D(\bar{\mathcal{E}}))$, which admits $\mu$ as an invariant measure.
\vskip.10in
From  Dirichlet form theory we obtain  easily:
\vskip.10in
\th{Corollary 1.4} $(\bar{\mathcal{E}}, D(\bar{\mathcal{E}}))$ and $({\mathcal{E}}, D({\mathcal{E}}))$ are recurrent in the sense that their associated semigroups $(T_t^i)_{t>0}, i=1,2,$ satisfy for $i=1,2$ $$\int_0^\infty T_t^ifdt =0 \textrm{ or } \infty \textrm{ a.e. for any } f\in L^1(E;\mu)\textrm{ with }f\geq0.$$
Here we use $(T_t^i)_{t>0}$ to denote the semigroup associated with the above Dirichlet forms respectively.
\vskip.10in
Recently a new uniform estimate for the solution $\Phi$ has been obtained in [MW17], which combined with the strong Feller property for $\Phi$ obtained in [HM16] and a support theorem in [HS17] for $\Phi$, may imply the exponential convergence to equilibrium in this case. By this result we can deduce the following estimate by using Dirichlet form constructed above. 

\vskip.10in
\th{Corollary 1.5} Suppose that the exponential convergence in the $L^2$-sense hold for the semigroup $\bar{P}_t$ associated with the solution $\Phi$. 
Then the following Poincar\'{e} inequality holds:
$$\mu(f^2)\leq C\mathcal{E}(f,f)+\mu(f)^2, \quad f\in D(\mathcal{E})$$
for some $C>0$. Moreover, there exists $c_0>0$ such that 
$$\int e^{c_0\|x\|_E}\mu(dx)<\infty,$$
where $E$ is the state space we introduced in Section 4. 
\vskip.10in

\th{Remark 1.6} In fact, Poincar\'{e} inequality implies the irreducibility of the Dirichlet form $(\mathcal{E},D(\mathcal{E}))$. Then by Corollary 1.4 and [FOT94, Theorem 4.7.1], for any nearly Borel non-exceptional set $B$,
$$P^x(\sigma_B\circ\theta_n<\infty,\forall n\geq0)=1, \quad \textrm{ for q.e. } x\in E.$$
Here $\sigma_B=\inf\{t>0:\Phi_t\in B\}$, $\theta$ is the shift operator for the Markov process $\Phi$, and for the definition of any nearly Borel non-exceptional set we refer to [FOT94]. Moreover by [FOT94, Theorem 4.7.3] we obtain the following strong law of large numbers: for $f\in L^1(E,\mu)$
$$\lim_{t\rightarrow\infty}\frac{1}{t}\int_0^tf(\Phi_s)ds=
\int fd\mu, \quad P^{x}-a.s.,$$
for q.e. $x\in E$.
\vskip.10in
\th{Remark 1.7} From Theorem 1.3 we know that there exists another Markov process which admits $\mu$ as an invariant measure. Is this Markov process the same as the solution $\Phi$ to (1.2) obtained in [MW16]? In Dirichlet form theory it corresponds to the problem of the relations between the domains of the Dirichlet forms $D(\mathcal{E})$ and $D(\bar{\mathcal{E}})$.
In the two dimensional case, they are the same (corresponding to restricted Markov uniqueness, see [RZZ15]). In the three dimensional case we do not know the answer until now, since the measure is more singular and  we do not know along which vector fields the integration by parts formula holds. This is also a major problem in Dirichlet form theory, which is  related to the long-standing open problem whether Markov uniqueness holds for the associated generator.

\vskip.10in

 The structure of this paper is as follows. In Section 2 we prove some useful estimates for the solutions to (1.2). In Section 3 we recall the lattice approximations, which is required to prove $\Phi$ is reversible w.r.t. $\mu$.   In Section 4 we give the proof of our first main result. In Section 5 we identify the Dirichlet form on the cylinder functions. In Appendix A, we recall some basic notions and results for the paracontrolled distribution method. In Appendix B, we calculate the convergence of the stochastic terms. We recall the paracontrolled analysis for the solutions to the lattice approximations in Appendix C. We also recall the definitions of Markov processes and quasi-regular Dirichlet forms  in Appendix D.

\vskip.20in
\no Notations: Let $\mathcal{S}'(\mathbb{T}^d)$ be the space of distributions on $\mathbb{T}^d=[-1,1]^d$. For $\alpha\in \mathbb{R}$, the H\"{o}lder-Besov space $\mathcal{C}^\alpha$ is given by $\mathcal{C}^\alpha=B^\alpha_{\infty,\infty}(\mathbb{T}^d)$ and for $p>1$ we use the notation $B^\alpha_p:=B^\alpha_{p,\infty}$.
For the definition of the general Besov spaces $B^\alpha_{p,q}$ and the paraproduct see  Appendix A. For $\beta>0,\alpha\in\mathbb{R}$ we write $\|\cdot\|_{\alpha}$, $C_T\mathcal{C}^\alpha$ and $C_T^\beta\mathcal{C}^\alpha$ instead of $\|\cdot\|_{B^\alpha_{\infty,\infty}}$,  $C([0,T];\mathcal{C}^\alpha)$ and $C^\beta([0,T];\mathcal{C}^\alpha)$, respectively in the following for simplicity. For a Banach space $E$, $\mathcal{B}(E)$ denotes the Borel-algebra on $E$ and $C_b(E)$ and $\mathcal{B}_b(E)$ denote the bounded continuous function and the bounded measurable functions on $E$, respectively.
The Fourier transform and the inverse Fourier transform are denoted by $\mathcal{F}$ and $\mathcal{F}^{-1}$. The heat semigroup is denoted by $P_t:=e^{t\Delta}$.

For $f\in \mathcal{S}'(\mathbb{T}^3)$ we write $\rho_\varepsilon*f:=\sum_kg(\varepsilon k)\langle f,e_k\rangle e_k$ with $g$ being a smooth radical function with compact support and $g(0)=1$, $g(\varepsilon k)=\mathcal{F}\rho_\varepsilon(k)$. Here and in the following $\langle\cdot,\cdot\rangle$ denotes $L^2(\mathbb{T}^3)$-inner product and $e_k(\xi)=2^{-3/2}e^{\iota\pi k\cdot \xi}$ for $k=(k^1,k^2,k^3)\in \mathbb{Z}^3,\xi=(\xi^1,\xi^2,\xi^3)\in\mathbb{T}^3$.
We also use $|k|_\infty=\max(|k^1|,|k^2|,|k^3|)$ and $\delta_{st}f:=f(t)-f(s)$. To make our paper better readable we summarize the graph notation used in the paper in the following table. The definition of them will be introduced below.
\begin{table}[!htbp]
\begin{tabular}{|c|c|c|c|c|c|c|c|c|c|}
\hline
\hline
$\Phi_1$ &  $\bar{\Phi}_1^\varepsilon$ & $-\Phi_2$ & $-\bar{\Phi}_2^\varepsilon$ & $-\rho_\varepsilon*\Phi_2$ & $({\Phi}_1)^{\diamond,2}$ & $(\bar{\Phi}^\varepsilon_1)^{\diamond,2}$  \\
\hline
$\includegraphics[height=0.5cm]{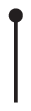}$ & $\includegraphics[height=0.5cm]{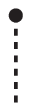}$ & $\includegraphics[height=0.7cm]{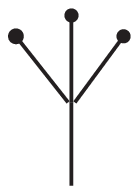}$ & $\includegraphics[height=0.7cm]{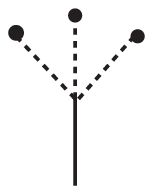}$ & $\includegraphics[height=0.7cm]{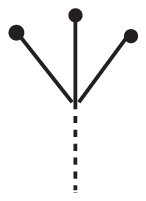}$ & $\includegraphics[height=0.5cm]{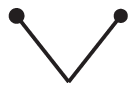}$ &$\includegraphics[height=0.5cm]{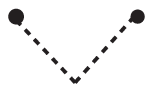}$   \\
\hline
 $K$ & $\bar{K}^\varepsilon$ & $\rho_\varepsilon*K$ & $(\rho_\varepsilon*\Phi_1)^{\diamond,3}$& $\Phi_1\diamond\Phi_2$ & $(\Phi_1)^{\diamond,2}\diamond\Phi_2$&  \\
\hline
  $\includegraphics[height=0.7cm]{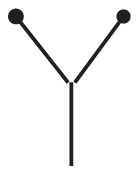}$ & $\includegraphics[height=0.7cm]{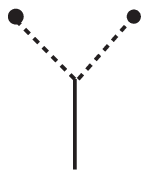}$& $\includegraphics[height=0.7cm]{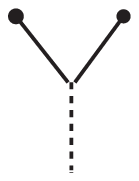}$& $\includegraphics[height=0.5cm]{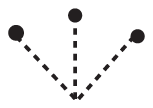}$& $-\includegraphics[height=0.7cm]{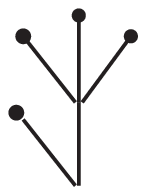} $&$-\includegraphics[height=0.7cm]{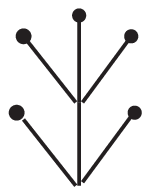} $&  \\
\hline
\end{tabular}
\end{table}

\section{A uniform estimate}
In this section we  give an uniform estimate of the solution to (1.2). In the following we assume that $\Phi_0\in \mathcal{C}^{-z}$ and $z\in(\frac{1}{2},\frac{2}{3})$.
 We fix $\kappa,\gamma>0$ satisfying
$$z-\frac{1}{2}>2\kappa,\quad 6\kappa<\gamma, \quad 10\kappa+3\gamma<2-3z.$$
 Parameters $\kappa,\gamma$ satisfying the above conditions can always be found. Indeed, we first choose $\gamma<\frac{2-3z}{3}$. Then the conditions are satisfied if we choose $\kappa>0$ small enough satisfying $\kappa<\frac{\gamma}{6}\wedge \frac{2z-1}{4}\wedge \frac{2-3z-3\gamma}{10}$.

Now we recall that the solution obtained by [CC13] and [MW16]: (1.2) can be split as follows:  $\Phi=\Phi_1+\Phi_2+\Phi_3$ and
$$ \Phi_1 (t)=\int_{-\infty}^t P_{t-s}  dW=\includegraphics[height=0.5cm]{diri01.eps},$$
$$\Phi_2 (t)=-\lim_{\varepsilon\rightarrow0}\int_0^tP_{t-s}\includegraphics[height=0.5cm]{diri11.eps}ds:=-\includegraphics[height=0.7cm]{diri03.eps},$$
and
$$\aligned \Phi_3 (t)=P_t(\Phi_0-\Phi_1 (0))-\int_0^t
P_{t-s}\bigg{[} & \Phi_3^3+3\Phi_3^2(\includegraphics[height=0.5cm]{diri01.eps}-\includegraphics[height=0.7cm]{diri03.eps})+\Phi_3(3(\includegraphics[height=0.7cm]{diri03.eps})^2
-6 \includegraphics[height=0.7cm]{diri14.eps})+3\includegraphics[height=0.5cm]{diri06.eps} \diamond\Phi_3\\& +3(\includegraphics[height=0.5cm]{diri01.eps} \diamond (\includegraphics[height=0.7cm]{diri03.eps} )^{2}-\includegraphics[height=0.7cm]{diri15.eps})-(\includegraphics[height=0.7cm]{diri03.eps} )^{3}  -(9\varphi-m)\Phi\bigg{]}ds.\endaligned\eqno(2.1)$$
Here we use $\includegraphics[height=0.5cm]{diri01.eps}$ to denote $\Phi_1$ and $\includegraphics[height=0.5cm]{diri02.eps}$ to denote $\rho_\varepsilon*\Phi_1$ and  introduce $\includegraphics[height=0.5cm]{diri06.eps}, \includegraphics[height=0.7cm]{diri03.eps}$ to reprensent $\Phi_1^{\diamond2}, -\Phi_2$, respectively.
$$\includegraphics[height=0.5cm]{diri06.eps},\includegraphics[height=0.5cm]{diri11.eps} , \includegraphics[height=0.7cm]{diri14.eps}, \includegraphics[height=0.5cm]{diri01.eps}\diamond (\includegraphics[height=0.7cm]{diri03.eps})^2,  \includegraphics[height=0.7cm]{diri15.eps},\varphi$$ involve a renormalization procedure and are defined in Appendix B. Throughout this paper we do not use the explicit formulation of these stochastic terms, but only use their regularity. We will introduce their regularity in (2.2) below.
 The most difficult part for renormalization is $\includegraphics[height=0.5cm]{diri06.eps}\diamond \Phi_3$. For this term we define $$K(t):=\int_0^tP_{t-s} (\Phi_1)^{\diamond,2}ds:=\includegraphics[height=0.7cm]{diri08.eps}.$$
 We have the following paracontrolled ansatz  $$\Phi_3 =-3\pi_<(-\includegraphics[height=0.7cm]{diri03.eps} +\Phi_3 ,\includegraphics[height=0.7cm]{diri08.eps} )+\Phi^{\sharp}$$
 with $\Phi^\sharp(t)\in \mathcal{C}^{1+3\kappa}$ for $t>0$. Here  $\Phi^\sharp$ is the regular term in the paracontrolled ansatz.
Then $$\aligned \includegraphics[height=0.5cm]{diri06.eps}\diamond \Phi_3:=&\pi_0(\Phi^{\sharp},\includegraphics[height=0.5cm]{diri06.eps})
-3C(-\includegraphics[height=0.7cm]{diri03.eps} +\Phi_3,\includegraphics[height=0.7cm]{diri08.eps} ,\includegraphics[height=0.5cm]{diri06.eps})
\\&-3(-\includegraphics[height=0.7cm]{diri03.eps} +\Phi_3)\pi_{0,\diamond}(\includegraphics[height=0.7cm]{diri08.eps},\includegraphics[height=0.5cm]{diri06.eps})+\pi_{<,>}(\Phi_3,\includegraphics[height=0.5cm]{diri06.eps}),\endaligned $$
where $C(-\includegraphics[height=0.7cm]{diri03.eps} +\Phi_3,\includegraphics[height=0.7cm]{diri08.eps} ,\includegraphics[height=0.5cm]{diri06.eps})$ is defined in Lemma A.3 and $\pi_{0,\diamond}(\includegraphics[height=0.7cm]{diri08.eps},\includegraphics[height=0.5cm]{diri06.eps})$ is defined in Appendix B.
 Now we introduce the following notations:$$\aligned C_W(T):=&\sup_{t\in[0,T]}\big[\|\includegraphics[height=0.5cm]{diri01.eps}\|_{-\frac{1}{2}-2\kappa}
 +\|\includegraphics[height=0.5cm]{diri06.eps}\|_{-1-2\kappa}+\|\includegraphics[height=0.7cm]{diri03.eps} \|_{\frac{1}{2}-2\kappa}+\|\pi_{0,\diamond}( \includegraphics[height=0.7cm]{diri03.eps} ,\includegraphics[height=0.5cm]{diri01.eps} )\|_{-2\kappa}\\&+\|\pi_{0,\diamond}( \includegraphics[height=0.7cm]{diri03.eps} ,\includegraphics[height=0.5cm]{diri06.eps})\|_{-\frac{1}{2}-2\kappa}+\|\pi_{0,\diamond}
(\includegraphics[height=0.7cm]{diri08.eps},\includegraphics[height=0.5cm]{diri06.eps})\|_{-2\kappa}\big]+\|\includegraphics[height=0.7cm]{diri03.eps} \|_{C^{\frac{1}{8}}_T\mathcal{C}^{\frac{1}{4}-2\kappa}},\endaligned\eqno(2.2)$$
and $$\rho_L:=\inf\{ t\geq0:C_W(t)\geq L\}.$$
 By [CC13] $P(C_W(T)<\infty, \forall T>0)=1$ and by [CC13] on this set there exists  a unique local solution $\Phi_3$ to (2.1). Recently in [MW16] the authors proved that the solution to (2.1) does not blow up in finite time. In fact we can check that the solution obtained in [MW16] satisfies (2.1) by smooth approximation.   In the following we consider the solution $\Phi$ obtained in [CC13] and [MW16].

\vskip.10in
Then we have the following  estimate for $\Phi$:
\vskip.10in
 \th{Proposition 2.1} For any $T>0$ there exist $C_0,\bar{m}>0$ depending on $L,T$ such that on the set $\{\rho_L>1\}$
 $$\sup_{t\in[0,T\wedge \rho_L]}[\|\Phi\|_{-z}+t^{\frac{\gamma+z+\kappa}{2}}\|\Phi_3\|_\gamma+t^{\frac{\frac{1}{2}+z+5\kappa}{2}}\|\Phi_3\|_{\frac{1}{2}+4\kappa}]\leq \exp{\{e^{C_0(\|\Phi_0\|_{-z}^{\bar{m}}+1)}\}}.$$

\vskip.10in
 \th{Remark} Here we obtain the estimate on the set $\{\rho_L>1\}$, since on this set we can choose $t^*$ below and the bound independent of $\omega$.
\vskip.10in
 \proof Set $$Q(t):=t^{\frac{\gamma+z+\kappa}{2}}\|\Phi_3\|_\gamma+t^{\frac{\frac{1}{2}+z+5\kappa}{2}}\|\Phi_3\|_{\frac{1}{2}+4\kappa}
 +t^{\frac{3(\gamma+z+\kappa)}{2}}\|\Phi^{\sharp}\|_{1+3\kappa}+1.$$
 By similar calculations as in [ZZ15, Section 4]  there exists $q>1$ such that for $t\leq \rho_L\wedge T$
 $$Q(t)^q\leq \bar{C}(\|\Phi_0\|_{-z}^q+1)+\bar{C}\int_0^tQ(s)^{3q}ds,$$
 where the constant $\bar{C}$ depends on $L,T,q$. Then Bihari's inequality implies that on the set $\{\rho_L>1\}$ for $t^*:=\bar{C}^{-1}[2 \bar{C}(\|\Phi_0\|_{-z}^q+1)]^{-2}\wedge 1\wedge T$
 $$\sup_{t\in [0,t^*]}Q(t)^q\leq C(\|\Phi_0\|_{-z}^q+1),$$
 where the constant $C$ depends on $L,T,q$. Then we obtain that
 $$\sup_{t\in [0,t^*]}[t^{\frac{\gamma+z+\kappa}{2}}\|\Phi_3\|_\gamma+t^{\frac{\frac{1}{2}+z+5\kappa}{2}}\|\Phi_3\|_{\frac{1}{2}+4\kappa}]\leq C(\|\Phi_0\|_{-z}+1).$$
 Moreover, by similar calculations as in [ZZ15, Section 4] there exists $m_0>0$ such that $$\sup_{t\in[0,t^*]}\|\Phi_3(t)\|_{-z}\leq  C(\|\Phi_0\|^{m_0}_{-z}+1),$$ and $$\|\Phi_3(t^*)\|_{\frac{1}{2}+4\kappa}\leq(t^*)^{-\frac{\frac{1}{2}+z+5\kappa}{2}}C(\|\Phi_0\|_{-z}+1)\leq C(\|\Phi_0\|^{m_0}_{-z}+1). $$
Consider the solution to (2.1) starting at $t^*$. By Proposition 2.2 we obtain that there exists some $m_1>0$ such that
$$\sup_{t\in [t^*,T\wedge \rho_L]}\|\Phi_3(t)\|_{\frac{1}{2}+4\kappa}\leq \exp{\{e^{C( \|\Phi_3(t^*)\|^{m_1}_{\frac{1}{2}+4\kappa}+1)}\}}.\eqno(2.3)$$
Thus the result follows.$\hfill\Box$

\vskip.10in
In the following proposition we use the result and notations from [MW16].
\vskip.10in
\th{Proposition 2.2} Let $\Phi_3$ be the solution to (2.1) with $\Phi_0-\Phi_1(0)$ replaced by $\Phi_3(0)\in \mathcal{C}^{\frac{1}{2}+4\kappa}$. Then there exists a constant $m_1>0$ such that on the set $\{\rho_L>1\}$  for any $T>0$
 $$\sup_{t\in[0,T\wedge \rho_L]}\|\Phi_3\|_{\frac{1}{2}+4\kappa}\leq \exp{\{e^{C(\|\Phi_3(0)\|_{\frac{1}{2}+4\kappa}^{m_1}+1)}\}}.$$
\vskip.10in
  Following  [MW16] we split the solution to (2.1) into the solutions to the following two equations:
 $$\left\{\begin{array}{ll}(\partial_t-\Delta)v=F(v+w)-cv,&\ \ \ \ v(0)=\Phi_3(0),\\(\partial_t-\Delta)w=G(v,w)+cv&\ \ \ \  w(0)=0,\end{array}\right.\eqno(2.4)$$
 with $$F(v+w):=-3\pi_<(v+w-\includegraphics[height=0.7cm]{diri03.eps},\includegraphics[height=0.5cm]{diri06.eps}),$$
 and $$G(v,w):=-(v+w)^3-3\textrm{com}(v,w)-3\pi_0(w,\includegraphics[height=0.5cm]{diri06.eps})-3\pi_>(v+w-\includegraphics[height=0.7cm]{diri03.eps},
 \includegraphics[height=0.5cm]{diri06.eps})+P(v+w),$$
 where $$\textrm{com}(v,w)=\pi_0(\textrm{com}_1(v,w),\includegraphics[height=0.5cm]{diri06.eps})
 +C(-3(v+w-\includegraphics[height=0.7cm]{diri03.eps}),\includegraphics[height=0.7cm]{diri08.eps},\includegraphics[height=0.5cm]{diri06.eps})$$ and $$P(v+w)=a_0+a_1(v+w)+a_2(v+w)^2,$$ with $$\textrm{com}_1(v,w)=P_tv(0)-3\int_0^tP_{t-s}\pi_<(v+w-\includegraphics[height=0.7cm]{diri03.eps},\includegraphics[height=0.5cm]{diri06.eps})ds+3\pi_<(v+w
 -\includegraphics[height=0.7cm]{diri03.eps},\includegraphics[height=0.7cm]{diri08.eps})$$ and $$a_0=-m(\includegraphics[height=0.5cm]{diri01.eps}-\includegraphics[height=0.7cm]{diri03.eps})+(\includegraphics[height=0.7cm]{diri03.eps})^3
 -3\includegraphics[height=0.5cm]{diri01.eps}\diamond(\includegraphics[height=0.7cm]{diri03.eps})^2
 +3\pi_{0,\diamond}(\includegraphics[height=0.7cm]{diri03.eps},\includegraphics[height=0.5cm]{diri06.eps})
 -9\includegraphics[height=0.7cm]{diri03.eps}\pi_{0,\diamond}(\includegraphics[height=0.7cm]{diri08.eps},\includegraphics[height=0.5cm]{diri06.eps}),$$
 and $$a_1=-m+6\includegraphics[height=0.7cm]{diri14.eps}
 -3(\includegraphics[height=0.7cm]{diri03.eps})^2+9\pi_{0,\diamond}(\includegraphics[height=0.7cm]{diri08.eps},\includegraphics[height=0.5cm]{diri06.eps}),\quad a_2=-3\includegraphics[height=0.5cm]{diri01.eps}+3\includegraphics[height=0.7cm]{diri03.eps}.$$
 By [CC13] we obtain that $a_0,a_1,a_2\in C_T\mathcal{C}^{-\frac{1}{2}-2\kappa}$. In (2.4) we omit $\varphi$ for simplicity. This term does not cause any problem since $\sup_{t\in[0,T]}t^\rho|\varphi(t)|<\infty$ for any $\rho>0$. Here we  emphasize that we consider (2.4) before $\rho_L$ and the constant $c$ in (2.4) only depends on $L$. We start by proving the following lemma:
 \vskip.10in
\th{Lemma 2.3} On the set $\{\rho_L>1\}$  for any $T>0$ we have that for $\beta_0=\frac{1}{2}+4\kappa, \gamma_0=\frac{5}{4}+4\kappa$ with $2\kappa$ being the same as $\varepsilon$ in [MW16]
 $$\sup_{t\in[0,T\wedge \rho_L]}[\|v\|_{B^{\beta_0}_6}+\|w\|_{B^{\gamma_0}_2}]\leq e^{C(\|\Phi_3(0)\|_{\frac{1}{2}+4\kappa}^{m_1}+1)}.\eqno(2.5)$$

\proof We would like to obtain how explicitly the solutions $(v,w)$ depend on the initial value $\Phi_3(0)$. The estimates in [MW16, Sections 3-5] depend polynomially  on the initial condition. Thus, we explicitly calculate in the following how the estimates in Sections 6 and 7 depend on the initial condition. Following the proof of Theorem 6.1 in [MW16], we first prove that on the set $\{\rho_L>1\}$ there exists some $m_2>0$ such that for some $t_*>0$
$$\int_0^{t_*}\|w(r)\|^2_{B^{1+4\kappa}_2}dr+\int_0^{t_*}\|w(r)\|_{L^6}^6dr+\int_0^{t_*}\|v(r)\|_{B^{\beta_0}_6}^6dr\leq C(\|\Phi_3(0)\|_{\frac{1}{2}+4\kappa}^{m_2}+1).\eqno(2.6)$$
In the following  the constants we omit in writing $\lesssim$ do not depend on the initial value. By (2.4) and Lemmas A.2- A.4 and a similar calculation as in the proof of [MW16, Lemma 2.3] we have for $t\in[0,T\wedge \rho_L]$
$$\|v(t)\|_{B^{\beta_0}_6}\lesssim \|\Phi_3(0)\|_{\frac{1}{2}+4\kappa}+\int_0^t\frac{1}{(t-s)^{
\frac{\beta_0+1+2\kappa}{2}}}(\|v(s)\|_{B^{\beta_0}_6}+\|w(s)\|_{B^{\gamma_0}_2}+1)ds,\eqno(2.7)$$
and
$$\aligned\|w(t)\|_{B^{\gamma_0}_2}\lesssim& \int_0^t\frac{1}{(t-s)^{\frac{\gamma_0}{2}+\frac{1}{4}+2\kappa}}(\|v(s)\|^3_{B^{\beta_0}_6}+\|w(s)\|^3_{B^{\gamma_0}_2}+1)ds
\\&+\int_0^t\frac{1}{(t-s)^{\frac{\gamma_0}{2}}}(\|v(s)\|_{B^{\beta_0}_6}+\|w(s)\|_{B^{\gamma_0}_2}+1+s^{-\frac{1}{4}}\|\Phi_3(0)\|_{\frac{1}{2}+4\kappa})ds
\\&+\int_0^t\frac{1}{(t-s)^{\frac{\gamma_0}{2}}}\int_0^s\frac{1}{(s-r)^{\frac{3}{4}+2\kappa}}(\|v(r)\|_{B^{\beta_0}_6}+\|w(r)\|_{B^{\gamma_0}_2})drds
\\&+\int_0^t\frac{1}{(t-s)^{\frac{\gamma_0}{2}}}\int_0^s\frac{1}{(s-r)^{1+4\kappa}}[\|\delta_{rs}v\|_{L^2}+\|\delta_{rs}w\|_{L^2}]drds.\endaligned\eqno(2.8)$$
Here and in the following the constants we omit depend on $L, T$.
By changing the  order of the integrals the third term in (2.8) equals to the following:
$$C\int_0^t(t-r)^{\frac{1}{4}-2\kappa-\frac{\gamma_0}{2}}(\|v(r)\|_{B^{\beta_0}_6}+\|w(r)\|_{B^{\gamma_0}_2})dr$$
By [MW16, Theorem 3.1] and H\"{o}lder's inequality the term containing $v$ in the last line of (2.8) is bounded by a constant times the following:
$$\aligned &\int_0^t\frac{1}{(t-s)^{\frac{\gamma_0}{2}}}\int_0^s\frac{1}{(s-r)^{\frac{5}{6}+6\kappa}}[\|v(r)\|_{B^{\beta_0}_2}+1+(\int_0^t\|w(u)\|_{L^2}^3du)^{1/3}dr]ds
\\\lesssim&1+\int_0^t(t-r)^{\frac{1}{6}-6\kappa-\frac{\gamma_0}{2}}\|v(r)\|_{B^{\beta_0}_2}dr+(\int_0^t\|w(u)\|_{B^{\gamma_0}_2}^3du)^{1/3},\endaligned$$
where in the last step we  change the order of integrals.
Since Lemma A.1 implies that $\|\cdot\|_{L^6}\lesssim \|\cdot\|_{B^{\gamma_0}_2}$, by [MW16, Theorem 4.1] and H\"{o}lder's inequality the term containing $w$ in the last line of (2.8) is bounded by a constant times the following:
$$\aligned&\int_0^t\frac{1}{(t-s)^{\frac{\gamma_0}{2}}}\int_0^s\frac{1}{(s-r)^{\frac{7}{8}+4\kappa}}[1+\|\Phi_3(0)\|^3_{\frac{1}{2}+4\kappa}
+\|w(r)\|_{B^{\gamma_0}_2}]drds
\\&+1+(\int_0^t\|w(u)\|_{B^{\gamma_0}_2}^6du)^{\frac{1}{2}}+(\int_0^t\|w^2(u)\|_{B^1_2}^2du)^{\frac{1}{2}}
\\\lesssim &1+\|\Phi_3(0)\|^3_{\frac{1}{2}+4\kappa}+\int_0^t(t-r)^{\frac{1}{8}-\frac{\gamma_0}{2}-4\kappa}\|w(r)\|_{B^{\gamma_0}_2}dr\\&+(\int_0^t\|w(u)\|_{B^{\gamma_0}_2}^6du)^{\frac{1}{2}}+(\int_0^t\|w^2(u)\|_{B^1_2}^2du)^{\frac{1}{2}}.\endaligned\eqno(2.9)$$
Here in the inequality we  change the order of integrals. For the last term of (2.9) by [MW16, (5.31)] and [MW16, Theorem 5.1] we have
$$\aligned (\int_0^t\|w^2(u)\|_{B^1_2}^2du)^{\frac{1}{2}}\lesssim& (\int_0^t[\|w^2(s)\|_{L^2}^2+\|w\nabla w(s)\|_{L^2}^2]ds)^{\frac{1}{2}}
\\\lesssim& 1+\|\Phi_3(0)\|^3_{\frac{1}{2}+4\kappa}+\int_0^t\|w(s)\|_{B^{\gamma_0}_2}^2ds.\endaligned$$
Now set $L_1(t):=\|v(t)\|_{B^{\beta_0}_6}+\|w(t)\|_{B^{\gamma_0}_2}+1$. Then by the calculations above and H\"{o}lder's inequality we obtain that there exists some $q>2$ such that for $t\in[0,T\wedge \rho_L]$
$$L_1(t)^q\leq \bar{C}(\|\Phi_3(0)\|_{\frac{1}{2}+4\kappa}^{3q}+1)+\bar{C}\int_0^tL_1(s)^{3q}ds.$$
Here the constant $\bar{C}$ depends on $L,T,q$. Thus Bihari's inequality implies that  on the set $\{\rho_L>1\}$ for $t\leq t^*:=\bar{C}^{-1}[2 \bar{C}(\|\Phi_3(0)\|_{\frac{1}{2}+4\kappa}^{3q}+1)]^{-2}\wedge 1\wedge T$
 $$\sup_{t\in [0,t^*]}L_1(t)^q\leq C(\|\Phi_3(0)\|_{\frac{1}{2}+4\kappa}^{3q}+1).$$
Here the constant $C$ depends on $L,T,q$. Now by taking $t_*$ satisfying [MW16, Proposition 6.2] and being smaller than $t^*$, we obtain (2.6), since by Lemma A.1 $\|\cdot\|_{L^6}\lesssim \|\cdot\|_{B^{\gamma_0}_2}$.  Then by the proof of Theorem 6.1 in [MW16] we obtain that on the set $\{\rho_L>1\}$ for $(k+2)t_*\leq\rho_L\wedge T$ with $k\in\mathbb{N}$
 $$\aligned & \int_{(k+1)t_*}^{(k+2)t_*}\|w(r)\|^2_{B^{1+4\kappa}_2}dr+\int_{(k+1)t_*}^{(k+2)t_*}\|w(r)\|_{L^6}^6dr+\int_{(k+1)t_*}^{(k+2)t_*}\|v(r)\|_{B^{\beta_0}_6}^6dr\\\leq &C+\frac{1}{t_*}\bigg(\int_{kt_*}^{(k+1)t_*}\|w(r)\|^2_{B^{1+4\kappa}_2}dr+\int_{kt_*}^{(k+1)t_*}\|w(r)\|_{L^6}^6dr+\int_{kt_*}^{(k+1)t_*}
 \|v(r)\|_{B^{\beta_0}_6}^6dr\bigg)\\\leq
 &[C(\|\Phi_3(0)\|_{\frac{1}{2}+4\kappa}^{3q}+1)]^{2k+3},\endaligned$$
 which implies that there exists $m_3>0$ such that
$$\int_0^{T\wedge \rho_L}\|w(r)\|^2_{B^{1+4\kappa}_2}dr+\int_0^{T\wedge \rho_L}\|w(r)\|_{L^6}^6dr+\int_0^{T\wedge \rho_L}\|v(r)\|_{B^{\beta_0}_6}^6dr\leq e^{C(\|\Phi_3(0)\|_{\frac{1}{2}+4\kappa}^{m_3}+1)}.$$
Thus, the results follow from the iteration arguments in [MW16, Section 7].$\hfill\Box$
\vskip.10in

\no \emph{Proof of Proposition 2.2} Since by Lemma A.1 $\|\cdot\|_{L^\infty}\lesssim\|\cdot\|_{\frac{1}{2}+4\kappa}\lesssim \|\cdot\|_{B^{\gamma_0}_4},$ and by [MW16, Theorem 3.1] $\sup_{t\in[0,T\wedge\rho_L]}\|v\|_{\frac{1}{2}+4\kappa}$ can be controlled by $C(1+\|\Phi_0\|_{\frac{1}{2}+4\kappa}+\sup_{t\in[0,T\wedge\rho_L]}\|w\|_{L^\infty})$, it is sufficient to prove that
$$\sup_{t\in[0,T\wedge \rho_L]}\|w(t)\|_{B^{\gamma_0}_4}\leq \exp{\{e^{C(\|\Phi_3(0)\|_{\frac{1}{2}+4\kappa}^{m_1}+1)}\}}.\eqno(2.10)$$
As in [MW16, Section 7] we write $w(t)=\sum_{j=1}^8\mathcal{W}_j(t)$ with
$$\aligned\mathcal{W}_1(t)+\mathcal{W}_2(t)=&-\int_0^tP_{t-s}(v+w)^3ds,
\\\mathcal{W}_3(t)=&-3\int_0^tP_{t-s}\pi_0(\textrm{com}_1(v,w),\includegraphics[height=0.5cm]{diri06.eps})ds,
\\\mathcal{W}_4(t)=&-3\int_0^tP_{t-s}\pi_0(w,\includegraphics[height=0.5cm]{diri06.eps})ds,
\\\mathcal{W}_5(t)=&\int_0^tP_{t-s}[a_2v^2]ds,
\\\mathcal{W}_6(t)=&2\int_0^tP_{t-s}[a_2vw]ds,
\\\mathcal{W}_7(t)=&\int_0^tP_{t-s}[a_2w^2]ds,
\\\mathcal{W}_8(t)=&\int_0^tP_{t-s}[...]ds,\endaligned$$
where $$...=-3C(-3(v+w-\includegraphics[height=0.7cm]{diri03.eps}),\includegraphics[height=0.7cm]{diri08.eps},\includegraphics[height=0.5cm]{diri06.eps})
-3\pi_>(v+w-\includegraphics[height=0.7cm]{diri03.eps},\includegraphics[height=0.5cm]{diri06.eps})+a_0+a_1(v+w)+cv.$$
Similarly as in [MW16, Section 7], we bound each term separately. For $t\in [0,T\wedge \rho_L]$ Lemma A.4 implies that
$$\|\mathcal{W}_1(t)+\mathcal{W}_2(t)\|_{B^{\gamma_0}_4}\lesssim \int_0^t\frac{1}{(t-r)^{\frac{\gamma_0}{2}}}(\|w(r)\|_{L^{12}}^3+\|v(r)\|_{L^{12}}^3)dr \leq e^{C(\|\Phi_3(0)\|_{\frac{1}{2}+4\kappa}^{m_1}+1)},$$
where in the second inequality we used that by Lemma A.1 $\|\cdot\|_{L^{12}}\lesssim \|\cdot\|_{B^{\gamma_0}_2}, \|\cdot\|_{L^{12}}\lesssim \|\cdot\|_{B^{\beta_0}_6}$ and Lemma 2.3.  For $t\in [0,T\wedge \rho_L]$ by Lemma A.4 we have
$$\aligned\|\mathcal{W}_4(t)\|_{B^{\gamma_0}_4}\lesssim& \int_0^t\frac{1}{(t-r)^{\frac{\gamma_0}{2}}}\|\pi_0(w,\includegraphics[height=0.5cm]{diri06.eps})(r)\|_{L^{4}}dr \lesssim \int_0^t\frac{1}{(t-r)^{\frac{\gamma_0}{2}}}\|w(r)\|_{B_{4}^{\gamma_0}}dr
,\endaligned$$
where in the second inequality we used Lemma A.2. Lemmas A.2 and A.4 imply that for $t\in [0,T\wedge \rho_L]$
$$\aligned\|\mathcal{W}_5(t)\|_{B^{\gamma_0}_4}\lesssim& \int_0^t\frac{1}{(t-r)^{\frac{\gamma_0}{2}+\frac{1}{4}+4\kappa}}\|v^2(r)\|_{B_4^{\frac{1}{2}+4\kappa}}dr \\\lesssim& \int_0^t\frac{1}{(t-r)^{\frac{\gamma_0}{2}+\frac{1}{4}+4\kappa}}\|v(r)\|_{B_6^{\frac{1}{2}+4\kappa}}\|v(r)\|_{L^{12}}dr
\\\leq &e^{C(\|\Phi_3(0)\|_{\frac{1}{2}+4\kappa}^{m_1}+1)},\endaligned$$
where in the last inequality we used Lemma 2.3 and that $\|\cdot\|_{L^{12}}\lesssim\|\cdot\|_{B_6^{\beta_0}}$. Also by Lemmas A.2 and A.4 we have that for $t\in [0,T\wedge \rho_L]$
$$\aligned\|\mathcal{W}_6(t)\|_{B^{\gamma_0}_4}\lesssim& \int_0^t\frac{1}{(t-r)^{\frac{\gamma_0}{2}+\frac{1}{4}+4\kappa}}\|vw(r)\|_{B_4^{\frac{1}{2}+4\kappa}}dr \\\lesssim& \int_0^t\frac{1}{(t-r)^{\frac{\gamma_0}{2}+\frac{1}{4}+4\kappa}}\|v(r)\|_{B_{6}^{\frac{1}{2}+4\kappa}}\|w(r)\|_{B_{12}^{\frac{1}{2}+4\kappa}}dr
\\\leq &e^{C(\|\Phi_3(0)\|_{\frac{1}{2}+4\kappa}^{m_1}+1)}\int_0^t\frac{1}{(t-r)^{\frac{\gamma_0}{2}+\frac{1}{4}+4\kappa}}\|w(r)\|_{B_{4}^{\gamma_0}}dr.\endaligned$$
Here in the last inequality we used that $\|\cdot\|_{B^{\beta_0}_{12}}\lesssim \|\cdot\|_{B^{\gamma_0}_4}$ and Lemma 2.3. Also  Lemmas A.2 and A.4 imply that for $t\in [0,T\wedge \rho_L]$
$$\aligned\|\mathcal{W}_7(t)\|_{B^{\gamma_0}_4}\lesssim& \int_0^t\frac{1}{(t-r)^{\frac{\gamma_0}{2}+\frac{1}{4}+4\kappa}}\|w^2(r)\|_{B_4^{\frac{1}{2}+4\kappa}}dr \\\lesssim& \int_0^t\frac{1}{(t-r)^{\frac{\gamma_0}{2}+\frac{1}{4}+4\kappa}}\|w(r)\|_{B_{6}^{\beta_0}}\|w(r)\|_{L^{12}}dr
\\\leq &e^{C(\|\Phi_3(0)\|_{\frac{1}{2}+4\kappa}^{m_1}+1)},\endaligned$$
where in the last inequality we used that $\|\cdot\|_{L^{12}}\lesssim \|\cdot\|_{B^{\beta_0}_6} $ and Lemma 2.3. Again Lemmas A.2-A.4 imply that for $t\in [0,T\wedge \rho_L]$
$$\aligned\|\mathcal{W}_8(t)\|_{B^{\gamma_0}_4}\lesssim& \int_0^t\frac{1}{(t-r)^{\frac{\gamma_0}{2}+\frac{1}{4}+4\kappa}}\|...\|_{B_4^{-\frac{1}{2}-4\kappa}}dr \\\lesssim& \int_0^t\frac{1}{(t-r)^{\frac{\gamma_0}{2}+\frac{1}{4}+4\kappa}}(1+\|w(r)\|_{B_4^{\beta_0}}+\|v(r)\|_{B_{4}^{\beta_0}})dr
\\\leq &e^{C(\|\Phi_3(0)\|_{\frac{1}{2}+4\kappa}^{m_1}+1)},\endaligned$$
where in the last inequality we used that $\|\cdot\|_{B_4^{\beta_0}}\lesssim \|\cdot\|_{B_2^{\gamma_0}}$ and Lemma 2.3. For $\mathcal{W}_3$ we need more calculations: for $t\in [0,T\wedge \rho_L]$
$$\aligned\|\mathcal{W}_3(t)\|_{B^{\gamma_0}_4}\lesssim &\int_0^t\frac{1}{(t-r)^{\frac{\gamma_0}{2}}}\|\pi_0(\textrm{com}_1(v,w),\includegraphics[height=0.5cm]{diri06.eps})(r)\|_{L^4}dr \\\lesssim &\int_0^t\frac{1}{(t-r)^{\frac{\gamma_0}{2}}}\|\textrm{com}_1(v,w)(r)\|_{B^{1+4\kappa}_4}dr
\\\lesssim &1+\|\Phi_3(0)\|_{\frac{1}{2}+4\kappa}+\int_0^t\frac{1}{(t-r)^{\frac{\gamma_0}{2}}}\int_0^r\frac{1}{(r-s)^{1+4\kappa-\frac{\beta_0}{2}}}\|w\|_{B^{\beta_0}_4}dsdr
\\&+\int_0^t\frac{1}{(t-r)^{\frac{\gamma_0}{2}}}\int_0^r\frac{1}{(r-s)^{1+4\kappa}}\|\delta_{sr}w\|_{L^4}dsdr
\\\lesssim & e^{C(\|\Phi_3(0)\|_{\frac{1}{2}+4\kappa}^{m_1}+1)}+\int_0^t\frac{1}{(t-r)^{\frac{\gamma_0}{2}}}\int_0^r\frac{1}{(r-s)^{1+4\kappa}}
\|\delta_{sr}w\|_{L^4}dsdr,\endaligned$$
where in the third inequality we used [MW16, Lemma 4.3] and in the last inequality we used $\|\cdot\|_{B_4^{\beta_0}}\lesssim \|\cdot\|_{B_2^{\gamma_0}}$ and Lemma 2.3.
 Now we control $\|\delta_{sr}w\|_{L^4}$ by a similar calculation as in [MW16, Theorem 4.1]: Lemma A.6 implies that
$$\|\delta_{sr}w\|_{L^4}\lesssim (r-s)^{\frac{1}{8}}\|w(s)\|_{B^{\beta_0}_4}+\|\delta'_{sr}w\|_{L^4}\lesssim (r-s)^{\frac{1}{8}}e^{C(\|\Phi_3(0)\|_{\frac{1}{2}+4\kappa}^{m_1}+1)}+\|\delta'_{sr}w\|_{L^4},$$
where
 $$\delta'_{sr}w=w(r)-P_{r-s}w(s)=\sum_{j=1}^8(\mathcal{W}_j(r)-P_{r-s}\mathcal{W}_j(s)),$$
 and we used that $\|\cdot\|_{B^{\beta_0}_4}\lesssim \|\cdot\|_{B^{\gamma_0}_2}$ and Lemma 2.3 in the last inequality.
Then by Lemmas 4.2, 4.4 and 4.6 in [MW16] we obtain the estimate for $\sum_{j=1}^4\|(\mathcal{W}_j(r)-P_{r-s}\mathcal{W}_j(s))\|_{L^4}$. A similar calculation as above for $\mathcal{W}_j$ with $j=5,6,7,8$ implies that $\sum_{j=5}^8\|(\mathcal{W}_j(r)-P_{r-s}\mathcal{W}_j(s))\|_{L^4}$ is bounded. Combining all this we obtain that for $s,r\in [0,T\wedge \rho_L]$ $$\aligned\|\delta'_{sr}w\|_{L^4}\lesssim& (r-s)^{\frac{1}{8}}e^{C(\|\Phi_3(0)\|_{\frac{1}{2}+4\kappa}^{m_1}+1)}\big[1+(\int_0^r\|w\|_{B^{\gamma_0}_4}^{4}du)^{\frac{1}{4}}\big]
\\&+(r-s)^{\frac{1}{8}}e^{C(\|\Phi_3(0)\|_{\frac{1}{2}+4\kappa}^{m_1}+1)}\|w\|^{\frac{1}{2}}_{4,r},
\endaligned$$
where $\|w\|_{4,r}=\sup_{u,u'\leq r}\frac{\|\delta'_{u'u}w\|_{L^4}}{|u-u'|^{\frac{1}{8}}}.$ Then by using the fact that $x\leq a+b\sqrt{x}$ implies $x\lesssim a+b$ as in the proof of Theorem 4.1 of [MW16], we have
$$\aligned\frac{\|\delta'_{sr}w\|_{L^4}}{(r-s)^{\frac{1}{8}}}\lesssim&e^{C(\|\Phi_3(0)\|_{\frac{1}{2}+4\kappa}^{m_1}+1)}
\big[1+(\int_0^r\|w\|_{B^{\gamma_0}_4}^{4}du)^{\frac{1}{4}}\big].
\endaligned$$
Thus we obtain that  for $t\in [0,T\wedge \rho_L]$
$$\aligned\|\mathcal{W}_3(t)\|_{B^{\gamma_0}_4}\lesssim e^{C(\|\Phi_3(0)\|_{\frac{1}{2}+4\kappa}^{m_1}+1)}\big[1+(\int_0^t\|w\|_{B^{\gamma_0}_4}^{4}du)^{\frac{1}{4}}\big].\endaligned$$
Combining all the estimates for $\mathcal{W}_j(t)$ and using H\"{o}lder's inequality we obtain that there exists some $q\geq4$ for $t\in [0,T\wedge \rho_L]$ such that
$$\|w(t)\|_{B^{\gamma_0}_4}^{q}\leq e^{C(\|\Phi_3(0)\|_{\frac{1}{2}+4\kappa}^{m_1}+1)}(1+\int_0^t\|w\|_{B^{\gamma_0}_4}^{q}du),$$
which implies (2.10) by Gronwall's inequality.$\hfill\Box$
\vskip.10in

 \section{Lattice approximation }
  In this section we will recall the lattice approximation in [ZZ15] for later use.
For $N\geq1$, let $\Lambda^N=\{-N,-(N-1),...,N\}^3$. Set $\varepsilon=\frac{2}{2N+1}$. Every point $k\in\Lambda^N$ can be identified with $\xi=\varepsilon k\in\Lambda_\varepsilon=\{\xi=(\xi^1,\xi^2,\xi^3)\in\varepsilon \mathbb{Z}^3:-1<\xi^1,\xi^2,\xi^3<1\}.$
We view $\Lambda_\varepsilon$ as a discretisation of the continuous three-dimensional torus $\mathbb{T}^3$ identified with $[-1,1]^3$. Then for $n\geq1$ we set $L^{2n}(\Lambda^\varepsilon):=\{\|f\|_{L^{2n}(\Lambda^\varepsilon)}^{2n}:=\sum_{x\in\Lambda^\varepsilon}\varepsilon^3|f(x)|^{2n}<\infty\}.$
(1.1) can be approximated by the following  lattice $\Phi^4_3$-field measure $\mu^\varepsilon(dx)$:
$$N^{-1}_{\varepsilon}\Pi_{\xi\in\Lambda_\varepsilon}dx_\xi\exp\bigg(-\varepsilon\sum_{|\xi_1-\xi_2|=\varepsilon,\xi_1,\xi_2\in\Lambda_\varepsilon}
(x(\xi_1)-x(\xi_2))^2+(3C_0^{\varepsilon}-9C_1^\varepsilon-m)\sum_{\xi\in\Lambda_\varepsilon}\varepsilon^3 x^2(\xi)-\frac{1}{2}\sum_{\xi\in\Lambda_\varepsilon}\varepsilon^3
x^4(\xi)\bigg),$$
where $N_\varepsilon$ is a normalization constant and we choose $C_0^\varepsilon,C_1^\varepsilon$ as in [ZZ15, Section 1].
The following  stochastic PDEs on $\Lambda_\varepsilon$ are the stochastic quantizations associated with the lattice $\Phi^4_3$-field measure:
$$\aligned d \Phi^{\varepsilon}(t)=&(\Delta_\varepsilon \Phi^{\varepsilon}(t)-(\Phi^{\varepsilon})^3(t) +(3C_0^{\varepsilon}-9C_1^\varepsilon-m)
\Phi^{\varepsilon}(t))dt\\&+ dW_N(t)
\\\Phi^\varepsilon(0)=&\Phi^\varepsilon_0,\endaligned\eqno(3.1)$$
where we fix
a cylindrical Wiener process in (1.2)  on $L^2(\mathbb{T}^3)$  given by $\sum_k\beta_ke_k(\xi)$ for $\xi\in\mathbb{T}^3$ and restrict it to
$L^2(\Lambda_\varepsilon)$ as $W_N(\xi)=\sum_{|k|_\infty\leq N}\beta_ke_k(\xi)$ for $\xi\in\Lambda_\varepsilon$, which is also a cylindrical Wiener process on
$L^2(\Lambda_\varepsilon)$. Here $\{\beta_k\}$ is a family of independent Brownian motions on $(\Omega,\mathcal{F},P)$. Also we take $\Phi_0^\varepsilon$ independent of $W$.  For $\xi\in\Lambda_\varepsilon$ define
$$\Delta_\varepsilon f(\xi):=\varepsilon^{-2}\sum_{y\in\Lambda_\varepsilon, y\sim \xi}(f(y)-f(\xi)),$$
where the nearest neighbor relation $\xi\sim y$ is to be understood with periodic boundary conditions on $\Lambda_\varepsilon$.
For $\Phi_0^\varepsilon$ satisfying $E\|\Phi_0^\varepsilon\|_{L^2(\Lambda_\varepsilon)}^2<\infty$ by [PR07, Theorem 3.1.1] there exists a unique solution $\Phi^\varepsilon$ to (3.1).

Following [MW14/ZZ15] we define a suitable extension of functions  defined on $\Lambda_\varepsilon$ onto all of the torus $\mathbb{T}^3$ (which we identify with the interval $[-1,1]^3$) in the following way:
$$\aligned\textrm{Ext} Y(\xi):=&\frac{1}{2^3}\sum_{k\in\{-N,...,N\}^3}\sum_{y\in\Lambda_\varepsilon}\varepsilon^3e^{\imath\pi k\cdot(\xi-y)}Y(y).\endaligned\eqno(3.2)$$
Now we extend the solutions of (3.1) to all of $\mathbb{T}^3$.   Let $u^\varepsilon=\textrm{Ext}\Phi^\varepsilon$ for simplicity. We have the following equation:
$$u^\varepsilon(t)=P_t^\varepsilon \textrm{Ext}\Phi_0^\varepsilon-\int_0^tP_{t-s}^\varepsilon Q_N[(u^\varepsilon)^3-(3C_0^{\varepsilon}-9C_1^\varepsilon-m) u^{\varepsilon}]ds+\int_0^tP_{t-s}^\varepsilon P_NdW.\eqno(3.3)$$
 where  $P_t^\varepsilon=\textrm{Ext}e^{t\Delta_\varepsilon}$ and $Q_N u(x)=P_N u(x)+\Pi_N u(x)$ with  $$P_N=\mathcal{F}^{-1}1_{|k|_\infty\leq N}\mathcal{F},$$ and $\Pi_N$ is defined for $u$ satisfying supp$\mathcal{F}u\subset \{k:|k|_\infty\leq 3N\}$
 $$\aligned\Pi_Nu(x)=&\sum_{i_1,i_2,i_3\in\{-1,0,1\},\sum_{j=1}^3i_j^2\neq0} e_N^{i_1i_2i_3}\mathcal{F}^{-1}1_{k\in P^{i_1i_2i_3} }
 \mathcal{F}u(x)\\=&\sum_{i_1,i_2,i_3\in\{-1,0,1\},\sum_{j=1}^3i_j^2\neq0}P_N[e^{i_1i_2i_3}_N u]\endaligned$$ with $P^{i_1i_2i_3}=\{k:k^ji_j>N \textrm{ if } i_j=-1,1;|k^j|\leq N,
 \textrm{ if } i_j=0\}$ is a rectangular division of $\mathbb{Z}^3\backslash \{k\in\mathbb{Z}^3, |k|_\infty\leq N\}$,  $e^{i_1i_2i_3}_N(\xi)=\Pi_{j=1}^3e^{-\imath\pi(2N+1)i_j\xi^j}$.

As in [ZZ15] we split (3.3) into the following three equations:
$$ u_1^{\varepsilon}(t)=\int_{-\infty}^t P_{t-s}^\varepsilon P_N dW,$$
$$u_2^{\varepsilon}(t)=-\int_0^tP_{t-s}^\varepsilon Q_N[(u_1^{\varepsilon})^{\diamond,3}]ds$$
and
$$\aligned u_3^{\varepsilon}(t)=&P^\varepsilon_t(\textrm{Ext}\Phi_0^\varepsilon-u_1^{\varepsilon}(0))-\int_0^t
P_{t-s}^\varepsilon\bigg{[} Q_N[6u_1^{\varepsilon} u_2^{\varepsilon}u_3^{\varepsilon}+3u_1^{\varepsilon}(u_3^{\varepsilon})^2+3u_1^{\varepsilon} (u_2^{\varepsilon})^{2}+(u_2^{\varepsilon}+u_3^{\varepsilon})^3]\\&+P_N[3(u_1^{\varepsilon})^{\diamond,2} \diamond(u_2^{\varepsilon}+u_3^{\varepsilon})+3e_N^{i_1i_2i_3}(u_1^{\varepsilon})^{\diamond,2} \diamond(u_2^{\varepsilon}+u_3^{\varepsilon})-(9\varphi^\varepsilon-m) u^\varepsilon]\bigg{]}ds.\endaligned\eqno(3.4)$$
Here the terms containing $\diamond$ are defined as in [ZZ15, Section 4].
For (3.4) we can do paracontrolled analysis as in [ZZ15, Section 4] and define the corresponding regular term $u^{\varepsilon,\sharp}$ in the paracontrolled ansatz. Also  we define $$C^\varepsilon_W(T), E^\varepsilon_W(T),A_N(T),D_N(T),\delta C_W^{\varepsilon}(T)$$ similarly as the corresponding stochastic  terms in [ZZ15]. Here for the completeness of the paper we include the definition of all these terms in Appendix C. Now we introduce the following definition:
$$\rho^\varepsilon_L:=\inf\{ t\geq0:C^\varepsilon_W(t)+E^\varepsilon_W(t)+A_N(t)+D_N(t)\geq L\},\eqno(3.5)$$
and $$\tau^\varepsilon_{C_0}:=\inf\{ t\geq0:t^{\frac{\gamma+z+\kappa}{2}}\|u^\varepsilon_3\|_\gamma+t^{\frac{\frac{1}{2}+z+5\kappa}{2}}\|u_3^\varepsilon\|_{\frac{1}{2}+4\kappa}\geq \exp{\{e^{C_0(\|\Phi_0\|_{-z}^{\bar{m}}+1)}\}}+1\},\eqno(3.6)$$
with $C_0,\bar{m}$ obtained in Proposition 2.1.

Now we obtain the following  estimate for the lattice approximations:
\vskip.10in
 \th{Proposition 3.1} We have on the set $\{\rho_L>1\}$, that for any $T>0$ there exists $C_1>0$ such that
 $$\aligned &\sup_{t\in[0,T\wedge \rho_L\wedge\rho^\varepsilon_L\wedge\tau^\varepsilon_{C_0}]}[\|u^\varepsilon-\Phi\|_{-z}+t^{\frac{\gamma+z+\kappa}{2}}\|u^\varepsilon_3-\Phi_3\|_\gamma
 +t^{\frac{\frac{1}{2}+z+5\kappa}{2}}\|u_3^\varepsilon-\Phi_3\|_{\frac{1}{2}+4\kappa}]
 \\\leq &C_1(\varepsilon^{\frac{\kappa}{2}}+\delta C_W^{\varepsilon}(T)+ E^\varepsilon_W(T)+ A_N(T)+D_N(T)+\|\textrm{Ext}\Phi_0^\varepsilon-\Phi_0\|_{-z}) \exp{\big\{\exp{\{ e^{C_1(\|\Phi_0\|_{-z}^{\bar{m}}+1)}\}}\big\}},\endaligned$$
  where the constant $C_1$ depends on $L,T$.

\proof Let $$L^\varepsilon(t):=t^{\frac{\gamma+z+\kappa}{2}}\|u^{\varepsilon}_3-\Phi_3\|_\gamma
+t^{\frac{\frac{1}{2}+z+5\kappa}{2}}\|u^{\varepsilon}_3-\Phi_3\|_{\frac{1}{2}+4\kappa}+t^{\frac{3(\gamma+z+\kappa)}{2}}\|u^{\varepsilon,\sharp}
-\Phi^{\sharp}\|_{1+3\kappa}.$$
Since the nonlinear terms are given by polynomials, by similar calculations as in [ZZ15] and Proposition 2.2 we have that on the set $\{\rho_L>1\}$ there exists $q>1$ such that for $t\in[0,T\wedge \rho_L\wedge\rho^\varepsilon_L\wedge\tau^\varepsilon_{C_0}]$
 $$\aligned L^\varepsilon(t)^q\leq& \exp{\{e^{C(\|\Phi_0\|_{-z}^{\bar{m}}+1)}\}}(\varepsilon^{\kappa/2}+\delta C_W^{\varepsilon}(T)+ E^\varepsilon_W(T)+ A_N(T)+D_N(T)+\|\textrm{Ext}\Phi_0^\varepsilon-\Phi_0\|_{-z})^q\\&+\exp{\{e^{C(\|\Phi_0\|_{-z}^{\bar{m}}+1)}\}}\int_0^tL^\varepsilon(s)^{q}ds,\endaligned$$
 which by Gronwall's inequality implies that for $t\in[0,T\wedge \rho_L\wedge\rho^\varepsilon_L\wedge\tau^\varepsilon_{C_0}]$
 $$L^\varepsilon(t)\leq (\varepsilon^{\kappa/2}+\delta C_W^{\varepsilon}(T)+ E^\varepsilon_W(T)+ A_N(T)+D_N(T)+\|\textrm{Ext}\Phi_0^\varepsilon-\Phi_0\|_{-z})\exp{\big\{\exp{\{ e^{C(\|\Phi_0\|_{-z}^{\bar{m}}+1)}\}}\big\}},$$
 on $\{\rho_L>1\}$.
 Moreover, by similar calculations as in [ZZ15] we obtain that on $\{\rho_L>1\}$ for $t\in[0,T\wedge \rho_L\wedge\rho^\varepsilon_L\wedge\tau^\varepsilon_{C_0}]$
 $$\|u^\varepsilon(t)-\Phi(t)\|_{-z}\leq (\varepsilon^{\kappa/2}+\delta C_W^{\varepsilon}(T)+E^\varepsilon_W(T)+ A_N(T)+D_N(T)+\|\textrm{Ext}\Phi_0^\varepsilon-\Phi_0\|_{-z}) \exp{\big\{\exp{\{ e^{C(\|\Phi_0\|_{-z}^{\bar{m}}+1)}\}}\big\}}.$$
 $\hfill\Box$
\vskip.10in
Similarly as in the proof of [HM15, Corollary  1.2] we obtain the following estimate for the measure $\bar{\mu}^\varepsilon:=\mu^\varepsilon\circ \textrm{Ext}^{-1}$. Since $\mu^\varepsilon$ is a measure on $L^2(\Lambda^\varepsilon)$ and $\textrm{Ext}$ is an isometry from $L^2(\Lambda_\varepsilon)$ to $P_NL^2(\mathbb{T}^3)$, $\bar{\mu}^\varepsilon$ has full support on $P_NL^2(\mathbb{T}^3)$:
\vskip.10in
\th{Lemma 3.2} Let $n\in\mathbb{N}$. Then  there exists a constant $C$ independent of $\varepsilon$ such that
 $$\int\|x\|^{2n}_{-z}\bar{\mu}^\varepsilon(dx)\leq C.$$
Moreover, $\bar{\mu}^\varepsilon$ weakly converges to $\mu$ on $\mathcal{C}^{-z}$.

 \proof The following calculations on $\Lambda_\varepsilon$ essentially follow [MW14, Lemma 8.4]. Suppose $supp \theta\subset \{a\leq|k|\leq b\}$ for $\theta$ as in Appendix A and $a,b>0$. If $2^ja>\sqrt{3}N$, then $\int\|\Delta_jx\|_{L^{2n}(\mathbb{T}^3)}^{2n}\bar{\mu}^\varepsilon(dx)=0$. For $x\in\textrm{supp}\bar{\mu}^\varepsilon$ we have $$\Delta_jx=\sum_{|k|_\infty\leq N}\theta_j(k)\langle x,e_k\rangle e_k=\sum_{|k|_\infty\leq N}\theta_j(k)\langle \textrm{Ext}^{-1}x,e_k\rangle_\varepsilon e_k,$$
 where $\theta_j(\cdot):=\theta(2^{-j}\cdot)$ and $\langle\cdot,\cdot\rangle, \langle\cdot,\cdot\rangle_\varepsilon$ denote the inner products in  $L^2(\mathbb{T}^3)$ and $L^2(\Lambda_\varepsilon)$, respectively. Here we can take $\textrm{Ext}^{-1}$ since $\textrm{Ext}$ is an isometry from $L^2(\Lambda_\varepsilon)$ to $P_NL^2(\mathbb{T}^3)$. If $2^jb<N-1$,
then by changing variables we have
$$\aligned &\int\|\Delta_jx\|_{L^{2n}(\mathbb{T}^3)}^{2n}\bar{\mu}^\varepsilon(dx)\\=&\int \|\sum_{|k|_\infty\leq N}\theta_j(k)\langle x,e_k\rangle_\varepsilon e_k\|_{L^{2n}(\mathbb{T}^3)}^{2n}\mu^\varepsilon(dx)
\\=&2^{-3n}\int\sum_{y_i\in \Lambda_\varepsilon,i=1,...,2n}\varepsilon^{6n}\sum_{|k_i|_\infty\leq N,i=1,...,2n}\big(\Pi_{i=1}^{2n}\theta_j(k_i)e_{k_i}(\xi-y_i)\big)S^\varepsilon_{2n}(y_1,...,y_{2n})d\xi
\\=&C\int\sum_{y_i\in 2^j\Lambda_\varepsilon,i=1,...,2n}\varepsilon^{6n}2^{6nj}\sum_{|k_i|_\infty\leq 2^{-j}N,k_i\in2^{-j}\mathbb{Z}^3,i=1,...,2n}2^{-6jn}\big(\Pi_{i=1}^{2n}\theta(k_i)e^{\pi \iota k_i(2^j\xi-y_i)}\big)S^\varepsilon_{2n}(\frac{y_1}{2^j},...,\frac{y_{2n}}{2^j})d\xi
\\=&C\int \sum_{y_i\in 2^j\Lambda_\varepsilon,i=1,...,2n}\varepsilon^{6n}2^{6jn}\bigg(\Pi_{i=1}^{2n}\frac{1}{[1+2^{2j}\sum_{l=1}^32(1-\cos(\pi 2^{-j}(2^j\xi^l-y_i^l)))]^2}\bigg)S^\varepsilon_{2n}(\frac{y_1}{2^j},...,\frac{y_{2n}}{2^j})\\&\sum_{|k_i|_\infty\leq 2^{-j}N,k_i\in2^{-j}\mathbb{Z}^3,i=1,...,2n}2^{-6jn}\big(\Pi_{i=1}^{2n}\theta(k_i)(1-\underline{\Delta}_j)^2e^{\pi \iota k_i(2^j\xi-y_i)}\big)d\xi
\\=&C\int \sum_{y_i\in 2^j\Lambda_\varepsilon,i=1,...,2n}\varepsilon^{6n}2^{6jn}\bigg(\Pi_{i=1}^{2n}\frac{1}{[1+2^{2j}\sum_{l=1}^32(1-\cos(\pi 2^{-j}(2^j\xi^l-y_i^l)))]^2}\bigg)S^\varepsilon_{2n}(\frac{y_1}{2^j},...,\frac{y_{2n}}{2^j})\\&\sum_{|k_i|_\infty\leq 2^{-j}N,k_i\in2^{-j}\mathbb{Z}^3,i=1,...,2n}2^{-6jn}\big(\Pi_{i=1}^{2n}(1-\underline{\Delta}_j)^2\theta(k_i)e^{\pi \iota k_i(2^j\xi-y_i)}\big)d\xi
\\\leq &C\int( \sum_{y_1,y_2\in 2^j\Lambda_\varepsilon}\varepsilon^{6}2^{6j}\frac{1}{(1+|2^j\xi-y_1|^2)^2}\frac{1}{(1+|2^j\xi-y_2|^2)^2}(C^\varepsilon(\frac{y_1}{2^j},\frac{y_2}{2^j})+\lambda^2))^nd\xi
\\\lesssim& 2^{jn},\endaligned $$
where $S^\varepsilon_{2n}(y_1,...,y_{2n})$ is the $2n$ point function for $\mu^\varepsilon$ from [BFS83] and $C^\varepsilon$ is the covariance for the corresponding Gaussian measure on the lattice  and $$\underline{\Delta}_jf(k)=2^{2j}\sum_{k'\in2^{-j}\mathbb{Z}^3,k\sim k'}(f(k')-f(k)).$$ Here in the last equality we  use the integration by parts formula, since on the boundary $\theta$ vanishes and in the first inequality we used that the support of $\theta$ is contained in an annulus to count the number of non-zero terms and deduce
$$\big|\sum_{|k_i|_\infty\leq 2^{-j}N,k_i\in2^{-j}\mathbb{Z}^3,i=1,...,2n}2^{-6jn}\big(\Pi_{i=1}^{2n}(1-\underline{\Delta}_j)^2\theta(k_i)e^{\pi \iota k_i(2^j\xi-y_i)}\big)\big|\lesssim 1.$$
In addition, we use (8.2) and Theorem 6.1 in [BFS83] to control $S_{2n}^\varepsilon$ and the following:
when $\xi^1\in [-1,1], \frac{1}{1-\cos(\pi \xi^1)}\leq \frac{C}{(\xi^1)^2}$ and when $\xi^1\in [1,2], \frac{1}{1-\cos(\pi \xi^1)}=\frac{1}{1-\cos(\pi (\xi^1-2))}\leq \frac{C}{(\xi^1-2)^2}$ and when $\xi^1\in [-2,-1], \frac{1}{1-\cos(\pi \xi^1)}=\frac{1}{1-\cos(\pi (\xi^1+2))}\leq \frac{C}{(\xi^1+2)^2}$. Furthermore, in the last step we use that the covariance $C^\varepsilon(y_1,y_2)$ of the Gaussian measure is of order $|y_1-y_2|^{-1}$.

If $\frac{2^ja}{\sqrt{3}}\leq N\leq 2^jb+1$, we choose a smooth function $\chi$ which equals $1$ on $\{\frac{a}{2}\leq|k|\leq 4b\}$ and vanishes outside the annulus $\{\frac{a}{3}\leq|k|\leq 5b\}$. Let $\chi_j=\chi(2^{-j}\cdot)$. We have
$$\aligned &\int\|\Delta_jx\|_{L^{2n}(\mathbb{T}^3)}^{2n}\bar{\mu}^\varepsilon(dx)\\=&\int \|\sum_{k}\theta_j(k)\chi_j(k) \langle x,e_k\rangle e_k\|_{L^{2n}(\mathbb{T}^3)}^{2n}\bar{\mu}^\varepsilon(dx)\leq C\int \|\sum_{k}\chi_j(k)\langle x,e_k\rangle e_k\|_{L^{2n}(\mathbb{T}^3)}^{2n}\bar{\mu}^\varepsilon(dx)
\\\lesssim&N^{3}\int \|\sum_{k}\chi_j(k)\langle x,e_k\rangle_\varepsilon e_k\|_{L^{2n}(\Lambda_\varepsilon)}^{2n}\mu^\varepsilon(dx)
\\\lesssim&2^{3j}\sum_{\xi\in\Lambda_\varepsilon}\varepsilon^3\sum_{y_i\in \Lambda_\varepsilon,i=1,...,2n}\sum_{|k_i|_\infty\leq N,i=1,...,2n}\varepsilon^{6n}\big(\Pi_{i=1}^{2n}\chi_j(k_i)e_{k_i}(\xi-y_i)\big)S^\varepsilon_{2n}(y_1,...,y_{2n})\lesssim 2^{3j+jn}.\endaligned $$
Here in the second inequality we used Lemma C.2 and the estimate in the last inequality can be obtained by a similar argument as above and  the integration by parts formula holds for the periodic boundary conditions. Thus, the first result holds by choosing $n$ large enough and because of Lemma A.1. In fact, for any $\alpha<-\frac{1}{2},$ $\int\|x\|^{2n}_{\alpha}\bar{\mu}^\varepsilon(dx)\leq C.$ The second result follows from the tightness of the $\bar{\mu}^\varepsilon$ and from the fact that the corresponding Schwinger functions converge (see [P75] and [HM15, Corollary 1.2]). $\hfill\Box$

\section{Existence of the Dirichlet form}

Consider the normal filtration $(\mathcal{F}_t)_{t\geq0}$ generated by $W$. As we mentioned in Section 2, by [Hai14, CC13,  MW16] for every $x\in \mathcal{C}^{-z}$ there exists a unique solution $\Phi(x)$ to (1.2) starting from $x$.
By [HM15] we have that $\Phi$ satisfies the Markov property on $\mathcal{C}^{-z}$ with respect to the filtration $(\mathcal{F}_t)_{t\geq0}$. Define $$P^x(A):=P(\Phi(x)\in A).$$
$P^x$ is a measure on $\Omega':=C([0,\infty);\mathcal{C}^{-z})$ and we use $E^x$ to denote the expectation under $P^x$. We use $X$ to denote the canonical process on $\Omega'$ and equip $\Omega'$ by the natural filtration $(\mathcal{M}_t)_{t\geq0}$ generated by $X$ (cf. [MR92, Chapter IV, (1.7)]). We know $X$ has the same distribution as $\Phi$.
By the Markov property of $\Phi$ we know $(\Omega',\mathcal{M}:=\vee_{t\geq0}\mathcal{M}_t,(\mathcal{M}_t)_{t\geq0},X, P^x)_{x\in \mathcal{C}^{-z}}$ is also a Markov process (cf. Definition D.2). Here iii) in Definition D.2 follows from the measurablity of $x\mapsto\Phi(x)$.
 Now we prove the following:
\vskip.10in
\th{Lemma 4.1} $(\Omega',\mathcal{M},(\mathcal{M}_t)_{t\geq0},X, P^x)_{x\in \mathcal{C}^{-z}}$ is  a Feller process on $\mathcal{C}^{-z}$.

\proof It suffices to check that $E^xf(X(t))$ is a continuous function on $\mathcal{C}^{-z}$ for $f\in C_b(\mathcal{C}^{-z})$. We have
$$\aligned&|E^{x_1}f(X(t))-E^{x_2}f(X(t))|=|Ef(\Phi(t,x_1))-Ef(\Phi(t,x_2))|\\\leq& E|f(\Phi(t, x_1 ))-f(\Phi(t,{x_2}))|1_{t\leq \rho_L}+CP(t>\rho_L).\endaligned$$
Here $\Phi(x)$ denotes the solution to (1.2) starting from $x$ and $\rho_L$ is defined as in Section 2. The first term goes to zero as $x_1$ goes to $x_2$ in $\mathcal{C}^{-z}$ by [Hai14] and the second term goes to zero as $L$ goes to infinity since $EC_W(t)\leq C$ with $C_W$ defined in (2.2).
$\hfill\Box$
\vskip.10in
By  $P^x(X\in C([0,\infty);\mathcal{C}^{-z}))=1$ for $x\in \mathcal{C}^{-z}$ and by [Chung82, Section 2.3 Theorem 1] we know that the Feller process $(\Omega',\mathcal{M},(\mathcal{M}_t)_{t\geq0},X, P^x)_{x\in \mathcal{C}^{-z}}$  satisfies the corresponding strong Markov property (cf. iii) in Definition D.3).

 To construct the Dirichlet form associated with $X$, we first extend the Markov process to starting points from a larger space, which contains $L^2(\mathbb{T}^3)$ as a subspace.  Choose $E=H^{-z-\epsilon}:=B^{-z-\epsilon}_{2,2}$ with $\epsilon>0$ and $H=L^2(\mathbb{T}^3)$. By Lemma A.1 we have $\mathcal{C}^{-z}\subset E$  and the following relation holds:
$$E^{*}\subset H^*\backsimeq H\subset E.$$ In the following we use $\langle\cdot,\cdot\rangle,|\cdot|$ to denote the inner product and norm on $H$ respectively and $\langle\cdot,\cdot\rangle$ also denotes the dual relation between $E^*$ and $E$ if there is no confusion. Now we would like to extend $X$ to a process $X'$ with state space $E$ in such a way that each $x\in E\backslash \mathcal{C}^{-z}$ is a trap for $X'$
(see [MR92, page 118]). For notation's  simplicity we still use   $(\Omega',\mathcal{M},(\mathcal{M}_t)_{t\geq0}, X,P^x)_{x\in E}$ to denote $X'$. In the following $(\Omega',\mathcal{M},(\mathcal{M}_t)_{t\geq0}, X,P^x)_{x\in E}$ is a continuous strong Markov process with state space $E$.
Define the associated semigroup for $f\in \mathcal{B}_b(E), x\in E$ $$\bar{P}_tf(x):=E^xf(X(t)).$$ We also introduce the following cylinder functions
$$\mathcal{F}C_b^\infty=\{f_1(\langle l_1,\cdot\rangle,...,\langle l_m,\cdot\rangle)|m\in\mathbb{N}, f_1\in C_b^\infty(\mathbb{R}^m),l_1,...,l_m\in E^*\}.$$
Define for $f\in \mathcal{F}C_b^\infty$ and $l\in H$, $$\frac{\partial f}{\partial l}(z):=\frac{d}{ds}f(z+sl)|_{s=0},z\in E,$$
  that is, by the chain rule,
  $$\frac{\partial f}{\partial l}(z)=\sum_{j=1}^m\partial_jf_1(\langle l_1,z\rangle,\langle l_2,z\rangle,...,\langle l_m,z\rangle)\langle l_j,l\rangle_H.$$
Let $Df$ denote the $H$-derivative of $f\in \mathcal{F}C_b^\infty$, i.e. the map from $E$ to $H$ such that $$\langle Df(z),l\rangle=\frac{\partial f}{\partial l}(z)\textrm{ for all } l\in H, z\in E.$$
In the following we  prove that $\bar{P}_t$ is a symmetric semigroup with respect to $\mu$. For this we  use lattice approximation in Section 3 and  let $\Phi^\varepsilon(x)$ be the solution to (3.1) obtained in Section 3 starting from $x\in L^2(\Lambda_\varepsilon)$.
By  existence and uniqueness of the solutions to (3.1) and similar arguments as in [PR07, LR15, Section 4.3] we obtain that $\Phi^\varepsilon$ satisfies the Markov property w.r.t. $\{\mathcal{F}_t\}_{t\geq0}$. We define the semigroup of the lattice approximation: for $f\in C_b(L^2(\Lambda_\varepsilon)), x\in L^2(\Lambda_\varepsilon)$, $$\tilde{P}_t^\varepsilon f(x)=E(f(\Phi^\varepsilon(t,x))).$$ Since (3.1) is a gradient system, by [DZ02, Theorem 12.3.2] we have for $f,g\in C_b(L^2(\Lambda_\varepsilon))$
$$\int \tilde{P}_t^\varepsilon f(x)g(x)\mu^\varepsilon(dx)=\int f(x)\tilde{P}_t^\varepsilon g(x)\mu^\varepsilon(dx).\eqno(4.1)$$ We also define the semigroup for the extension of the lattice approximation on $P_NE$: for $f\in C_b(P_NE)$, $x\in P_NE$,
 $$\bar{P}_t^\varepsilon f(x)=E(f(u^\varepsilon(t,x))),$$
 where $P_N$ is as introduced in Section 3 and $u^\varepsilon(x)$ is the solution to (3.3) starting from $x$. Then we prove that $\bar{P}_t^\varepsilon$ is symmetric with respect to $\bar{\mu}^\varepsilon$. Since the extension operator Ext defined in (3.2) is an isometry from $L^2(\Lambda_\varepsilon)$ to $P_NE$,  we view $\bar{\mu}^\varepsilon$ as a measure on $P_NE$.
\vskip.10in
\th{Lemma 4.2} For $f, g\in \mathcal{F}C_b^\infty$ we have
$$\int \bar{P}_t^\varepsilon (f|_{P_NE})(x)g|_{P_NE}(x)\bar{\mu}^\varepsilon(dx)=\int f|_{P_NE}(x)\bar{P}_t^\varepsilon (g|_{P_NE})(x)\bar{\mu}^\varepsilon(dx),$$
where we used that $P_NE\subset E$.

\proof Without loss of generality we assume that $f(x)=f_1(\langle x,l\rangle), g(x)=g_1(\langle x,h\rangle)$ with $f_1,g_1\in C_b^\infty$.
Then we have that for $l_1=\sum_{|k|_\infty\leq N}\langle l,e_k\rangle e_k, h_1=\sum_{|k|_\infty\leq N}\langle h,e_k\rangle e_k,$
$$\aligned &\int \bar{P}_t^\varepsilon (f|_{P_NE})(x)g|_{P_NE}(x)\bar{\mu}^\varepsilon(dx)=\int E(f_1(\langle u^\varepsilon(t,x),l_1\rangle))g_1(\langle x,h_1\rangle)\bar{\mu}^\varepsilon(dx)\\=&\int E(f_1(\langle \Phi^\varepsilon(t,\textrm{Ext}^{-1}x),l_1\rangle_\varepsilon))g_1(\langle \textrm{Ext}^{-1}x,h_1\rangle_\varepsilon)\bar{\mu}^\varepsilon(dx)\\=&\int E(f_1(\langle \Phi^\varepsilon(t,x),l_1\rangle_\varepsilon))g_1(\langle x,h_1\rangle_\varepsilon)\mu^\varepsilon(dx)
\\=&\int E(g_1(\langle \Phi^\varepsilon(t,x),h_1\rangle_\varepsilon))f_1(\langle x,l_1\rangle_\varepsilon)\mu^\varepsilon(dx)=\int E(g_1(\langle u^\varepsilon(t,x),h_1\rangle))f_1(\langle x,l_1\rangle)\bar{\mu}^\varepsilon(dx)\\=&\int \bar{P}_t^\varepsilon (g|_{P_NE})(x)f|_{P_NE}(x)\bar{\mu}^\varepsilon(dx).\endaligned$$
Here in the second equality we used $\langle x,l_1\rangle=\langle \textrm{Ext}^{-1}x,l_1\rangle_\varepsilon$ for $x\in P_NE$ to deduce $\langle \Phi^\varepsilon_t,l_1\rangle_\varepsilon=\langle u^\varepsilon_t,l_1\rangle$ and in the forth equality we used (4.1).
$\hfill\Box$

\vskip.10in
By Lemma 4.2 and [MR92, Chapter II Prop. 4.3] we know that $(\bar{P}_t^\varepsilon)_{t>0}$ can be extended as a strongly continuous sub-Markovian semigroup of contractions on $L^2(P_NE;\bar{\mu}^\varepsilon)$. By [MR92, Chap I] there exists a corresponding Dirichlet form for $(\bar{P}_t^\varepsilon)_{t>0}$. In Proposition 4.4 we will give the explicit formula for this Dirichlet form.  Now we  prove that $\bar{P}_t$ is symmetric with respect to $\mu$.

\vskip.10in
\th{Proposition 4.3}For $f, g\in \mathcal{F}C_b^\infty$ we have for $t\geq0$
$$\int \bar{P}_t f(x)g(x)\mu(dx)=\int f(x)\bar{P}_tg(x)\mu(dx).$$

\proof By Lemma 4.2 it suffices to prove that for $f, g\in \mathcal{F}C_b^\infty$
$$\lim_{\varepsilon\rightarrow0}\int \bar{P}^\varepsilon_t (f|_{P_NE})(x)g|_{P_NE}(x)\bar{\mu}^\varepsilon(dx)=\int \bar{P}_t f(x)g(x)\mu(dx).\eqno(4.2)$$
Lemmas 3.2 and 4.1 imply that
$$\lim_{\varepsilon\rightarrow0}\int \bar{P}_t f(x)g(x)\bar{\mu}^\varepsilon(dx)=\int \bar{P}_t f(x)g(x)\mu(dx).$$
We also have
$$\aligned&\int |\bar{P}^\varepsilon_t (f|_{P_NE})(x)-\bar{P}_t f(x)||g(x)|\bar{\mu}^\varepsilon(dx)
\\\leq &C\int E\big(|f(u^\varepsilon(t,x))-f(\Phi(t,x))|1_{\{t<\rho_L^\varepsilon\wedge \rho_L, \rho_L>1\}}\big)\bar{\mu}^\varepsilon(dx)+C(P(t\geq\rho_L\wedge \rho_L^\varepsilon)+P(\rho_L\leq 1)),\endaligned\eqno(4.3)$$
where $\rho_L, \rho_L^\varepsilon$ are as introduced in Section 2 and (3.5), respectively.
The second term in (4.3) is bounded by a constant times
$$\aligned&P((C^\varepsilon_W+E^\varepsilon_W+A_N+D_N)(t)>L)+P(C_W(t)>L)+P(C_W(1)>L)
\\\leq &C/L,\endaligned$$
which uniformly goes to zero as $L$ goes to $\infty$.
For some $\delta_0>0$ the first term in (4.3) is bounded by
$$\aligned&\varepsilon^{\delta_0}C\int P(\|u^\varepsilon(t,x)-\Phi(t,x)\|_{-z} <\varepsilon^{\delta_0})\bar{\mu}^\varepsilon(dx)\\&+C\int P(t<\rho_L^\varepsilon\wedge \rho_L, \rho_L>1,\|u^\varepsilon(t,x)-\Phi(t,x)\|_{-z} >\varepsilon^{\delta_0})\bar{\mu}^\varepsilon(dx).\endaligned\eqno(4.4)$$
Then the first term is bounded by $C\varepsilon^{\delta_0}$ and the second integral in (4.4) is bounded by
$$\aligned &\int[P(t<\rho_L^\varepsilon\wedge \rho_L,\rho_L>1,t<\tau_{C_0}^\varepsilon, \|u^\varepsilon(t,x)-\Phi(t,x)\|_{-z} >\varepsilon^{\delta_0})+P(t<\rho_L^\varepsilon\wedge \rho_L,\rho_L>1,t\geq\tau_{C_0}^\varepsilon)]\bar{\mu}^\varepsilon(dx)\\\leq&2\int P(\sup_{s\in[0,\rho_L^\varepsilon\wedge \rho_L\wedge t\wedge\tau_{C_0}^\varepsilon]}[\|u^\varepsilon(x)-\Phi(x)\|_{-z} +s^{\frac{\gamma+z+\kappa}{2}}\|u^\varepsilon_3-\Phi_3\|_\gamma+s^{\frac{\frac{1}{2}+z+5\kappa}{2}}\|u_3^\varepsilon-\Phi_3\|_{\frac{1}{2}+4\kappa} ] >\varepsilon^{\delta_0},\\&\quad\quad\quad\rho_L>1)\bar{\mu}^\varepsilon(dx)
\\\leq& 2\int P(2C_1\varepsilon^{\kappa_0} \exp{\big\{\exp{\{ e^{C_1(\|x\|_{-z}^{\bar{m}}+1)}\}}\big\}}>\varepsilon^{\delta_0})\bar{\mu}^\varepsilon(dx)\\&+2\int P(\delta C_W^{\varepsilon}(t)+ A_N(t)+ E^\varepsilon_W(t)+D_N(t)>\varepsilon^{\kappa_0})\bar{\mu}^\varepsilon(dx)
\\\leq&2\int 1_{\{\|x\|_{-z}^{\bar{m}}>\frac{1}{C_1}\ln\ln\ln\frac{\varepsilon^{\delta_0-\kappa_0}}{2C_1}-1\}}\bar{\mu}^\varepsilon(dx)+2C\varepsilon^{\kappa_1-\kappa_0}\\\leq&2\int \frac{1}{\frac{1}{C_1}\ln\ln\ln\frac{\varepsilon^{\delta_0-\kappa_0}}{2C_1}-1}\|x\|_{-z}^{\bar{m}}\bar{\mu}^\varepsilon(dx)+2C\varepsilon^{\kappa_1-\kappa_0}\rightarrow0,\quad \textrm{ as } \varepsilon\rightarrow0,\endaligned$$
where $u^\varepsilon_3, \Phi_3$ correspond to $u^\varepsilon(x), \Phi(x)$ respectively and $\tau_{C_0}^\varepsilon$ is defined in (3.6) and in the first inequality we used Proposition 2.1 and the definition of $\tau_{C_0}^\varepsilon$ to deduce
 $$\sup_{s\in[0,\rho_L^\varepsilon\wedge \rho_L\wedge t\wedge\tau_{C_0}^\varepsilon]}[\|u^\varepsilon-\Phi\|_{-z} +s^{\frac{\gamma+z+\kappa}{2}}\|u^\varepsilon_3-\Phi_3\|_\gamma+s^{\frac{1/2+z+5\kappa}{2}}\|u_3^\varepsilon-\Phi_3\|_{\frac{1}{2}+4\kappa} ] >\varepsilon^{\delta_0}.$$
In the second inequality we used Proposition 3.1 and in the third inequality we used  Proposition C.1 and in the last step we used Lemma 3.2. Here we choose $0<\delta_0<\kappa_0<\kappa_1\wedge \frac{\kappa}{2}$ for $\frac{\kappa}{2}, \kappa_1$ coming from Proposition 3.1 and Proposition C.1, respectively. Summarizing, we obtain the result. $\hfill\Box$
\vskip.10in

Now we identify the Dirichlet form associated with $(\bar{P}_t^\varepsilon)_{t>0}$ on $L^2(P_NE,\bar{\mu}^\varepsilon)$.
\vskip.10in
\th{Proposition 4.4} The Dirichlet form associated with $(\bar{P}_t^\varepsilon)_{t>0}$ can be written as the closure of the following bilinear form
$$\mathcal{E}^\varepsilon(f,g)=\frac{1}{2}\sum_{|k|_\infty\leq N}\int_{P_NE}  \frac{\partial f}{\partial e_k} \frac{\partial g}{\partial e_k} d\bar{\mu}^\varepsilon, \quad f,g\in C_b^\infty(P_NE),$$
where $C_b^\infty(P_NE)$ means smooth functions on $P_NE$ with bounded derivatives.

\proof It is standard to obtain that the closure of $(\mathcal{E}^\varepsilon,C_b^\infty(P_NE))$ is a quasi-regular Dirichlet form (cf. Definition D.1, [MR92, Chap IV Section 4]), which is denoted by $(\mathcal{E}^\varepsilon,D(\mathcal{E}^\varepsilon))$. By Theorem D.4 there exists a Markov process with continuous sample paths properly associated with $(\mathcal{E}^\varepsilon,D(\mathcal{E}^\varepsilon))$.   Now we want to prove that the associated Markov process has the same distribution as $u^\varepsilon$.

We can easily conclude that the log-derivative of $\mu^\varepsilon$ along $e_k$ for $|k|_{\infty}\leq N$ is given by
$$b_k(x)=2\langle x,\Delta_\varepsilon e_k\rangle_\varepsilon -2\langle x^3-(3C_0^\varepsilon-9C_1^\varepsilon-m)x,e_k\rangle_\varepsilon\textrm{ for } x\in L^2(\Lambda_\varepsilon),$$
which implies that for $f\in C_b^\infty(P_NE)$ and $|k|_\infty\leq N$
$$\int \frac{\partial f}{\partial e_k}(x)d\bar{\mu}^\varepsilon=\int \frac{\partial }{\partial e_k}(f\circ \textrm{Ext})( x)d\mu^\varepsilon=-\int f(\textrm{Ext} x)b_k(x)d\mu^\varepsilon=-\int f(x)b_k(\textrm{Ext}^{-1}x)d\bar{\mu}^\varepsilon,$$
 we obtain that the log-derivative of $\bar{\mu}^\varepsilon$ is $$\beta_k(x)=b_k(\textrm{Ext}^{-1}x)=2\langle x,\Delta_\varepsilon e_k\rangle_{L^2(\mathbb{T}^3)}-2\langle Q_N(x^3-(3C_0^\varepsilon-9C_1^\varepsilon-m)x),e_k\rangle_{L^2(\mathbb{T}^3)},\quad x\in P_NE, |k|_\infty\leq N,$$
where we used that $\textrm{Ext}(\textrm{Ext}^{-1}x)^3=Q_N(x^3)$ for $x\in P_NE$.
 This implies that the associated Markov process is a probabilistically weak solution to the equation (3.3).  On the other hand, the equation (3.3) is a finite dimensional stochastic differential equation and we can easily obtain the pathwise uniqueness of the solutions to the equation (3.3). This deduces that $u^\varepsilon$ has the same distribution as the Markov process given by the Dirichlet form  $(\mathcal{E}^\varepsilon,D(\mathcal{E}^\varepsilon))$. By  Theorem D.4 we know that the semigroup of $u^\varepsilon$ $(\bar{P}_t^\varepsilon)_{t>0}$ is properly associated with $(\mathcal{E}^\varepsilon,D(\mathcal{E}^\varepsilon))$.
 $\hfill\Box$
\vskip.10in

\vskip.10in

\no\emph{Proof of Theorem 1.1}: By Proposition 4.3 we have that $\int \bar{P}_t fd\mu=\int f d\mu$ for $f\in \mathcal{F}C_b^\infty$. Since $\sigma(\mathcal{F}C_b^\infty)=\mathcal{B}(E)$, we deduce that $\mu$ is an invariant measure for the semigroup $\bar{P}_t$, which implies that
$$\int \bar{P}_t fd\mu=\int f d\mu \textrm{ for } f\in \mathcal{B}_b(E).\eqno(4.4)$$
By Proposition 4.3 and using (4.4) and the fact that $\mathcal{F}C_b^\infty$ is dense in $ L^2(E;\mu)$, we have that for $f, g\in \mathcal{B}_b(E)$
$$\int \bar{P}_t f(x)g(x)\mu(dx)=\int f(x)\bar{P}_tg(x)\mu(dx).$$
 Since $(\bar{P}_t)_{t>0}$ is  sub-Markovian, by [MR92, Chapter II Proposition 4.1] it can be extended to $L^2(E,\mu)$. This extension is still denoted by $(\bar{P}_t)_{t>0}$. On the other hand, since $\Phi$ has continuous path in $E$, we can deduce that $\bar{P}_tf \rightarrow_{t\rightarrow0} f$ in $\mu$-measure for $f\in \mathcal{F}C_b^\infty$. Then by [MR92, Chapter II Proposition 4.3] $(\bar{P}_t)_{t>0}$ is a strongly continuous contraction semigroup on $ L^2(E;\mu)$. Then there exists a corresponding Dirichlet form $(\mathcal{E}, D(\mathcal{E}))$ associated with $(\bar{P}_t)_{t>0}$.

We  know that $(\Omega', \mathcal{M}, (\mathcal{M}_t)_{t>0}, X, P^z )_{z\in E}$ is a right process in the sense of Definition D.3, which implies that  $(\mathcal{E}, D(\mathcal{E}))$ is a quasi-regular Dirichlet form by  Theorem D.4.

 In the following we prove that $\mathcal{F}C_b^\infty\subset D(\mathcal{E})$.
By (4.2) and since $\bar{\mu}^\varepsilon$ converges weakly to $\mu$ we know that for $f\in \mathcal{F}C_b^\infty$,
$$\aligned\sup_{t>0}\frac{1}{t}\int(\bar{P}_tf-f)fd\mu=&\sup_{t>0}\lim_{\varepsilon\rightarrow0}\frac{1}{t}\int(\bar{P}^\varepsilon_t(f|_{P_NE})-f|_{P_NE})f|_{P_NE}d\bar{\mu}^\varepsilon
\\\leq&\liminf_{\varepsilon\rightarrow0}\sup_{t>0}\frac{1}{t}\int(\bar{P}^\varepsilon_t(f|_{P_NE})-f|_{P_NE})f|_{P_NE}d\bar{\mu}^\varepsilon
\\=&\liminf_{\varepsilon\rightarrow0}\mathcal{E}^\varepsilon(f|_{P_NE},f|_{P_NE})<\infty,
\endaligned$$
where in the last inequality we used Proposition 4.4.
This implies that $\mathcal{F}C_b^\infty\subset D(\mathcal{E})$ and for $f\in \mathcal{F}C_b^\infty$,
$$\mathcal{E}(f,f)\leq\frac{1}{2} \int |Df|^2d\mu.\eqno(4.5)$$
For $l\in E^*$ by (4.5) we can easily find $f_n\in \mathcal{F}C_b^\infty$ such that $f_n\rightarrow \langle l,\cdot\rangle$ in $L^2(E,\mu)$ and $f_n$ is a Cauchy sequence in $D(\mathcal{E})$, which implies $\langle l,\cdot\rangle\in D(\mathcal{E})$ since $(\mathcal{E}, D(\mathcal{E}))$ is a closed form.
$\hfill\Box$

\section{Identification of the Dirichlet form}
In this section we identify the Dirichlet form $(\mathcal{E},D(\mathcal{E}))$ on $\mathcal{F}C_b^\infty$. To complete this, we first try to write the nonlinear term as an additive functional of the solution. Here we use paracontrolled analysis to prove the solution $\Phi$ to (1.2) satisfies the following equation in the analytic weak sense:
$$\Phi(t)= \Phi_0+\int_0^t\Delta \Phi ds-\lim_{\varepsilon\rightarrow0}\int_0^t[(\rho_\varepsilon\ast\Phi)^3-(3\bar{C}_0^\varepsilon\rho_\varepsilon*\Phi-9\tilde{C}_{1}^\varepsilon\Phi-m\Phi) ]ds +W(t),\eqno(5.1)$$
where $\bar{C}_0^\varepsilon$ and $\tilde{C}_1^\varepsilon$ are defined below.
For this we consider the following approximation: Let $\bar{\Phi}^\varepsilon$ be the solutions to the following equation:
 $$d \bar{\Phi}^{\varepsilon}=\Delta \bar{\Phi}^{\varepsilon}dt+\rho_\varepsilon* dW-(\bar{\Phi}^{\varepsilon})^3dt+(3 \bar{C}_0^\varepsilon-9\bar{C}_1^\varepsilon-m)\bar{\Phi}^{\varepsilon}dt,\eqno(5.2)$$
$$\bar{\Phi}^\varepsilon(0)=\Phi_0.$$
Here $\bar{C}_0^\varepsilon$ and $\bar{C}_1^\varepsilon$ are the corresponding constants defined in Appendix B. For this equation we can also write $\bar{\Phi}^\varepsilon=\bar{\Phi}^{\varepsilon}_1+ \bar{\Phi}^{\varepsilon}_2+\bar{\Phi}_3^\varepsilon$ and define $\bar{\Phi}^{\varepsilon}_1, \bar{\Phi}^{\varepsilon}_2$, $\bar{\Phi}_3^\varepsilon$, $\bar{K}^\varepsilon$, $\bar{\Phi}^{\varepsilon,\sharp}$ similarly as  in Section 2. Here we also introduce graph notations for them. We use $\includegraphics[height=0.5cm]{diri02.eps}$ to denote $\bar{\Phi}_1^\varepsilon$ and $\includegraphics[height=0.7cm]{diri04.eps}$ to denote $-\bar{\Phi}_2^\varepsilon$. Moreover, $\includegraphics[height=0.7cm]{diri09.eps}$ is used to denote $\bar{K}^\varepsilon$.  The corresponding renormalized terms $\includegraphics[height=0.5cm]{diri07.eps}$, $\includegraphics[height=0.5cm]{diri11.eps},\pi_{0,\diamond}(\includegraphics[height=0.7cm]{diri04.eps},\includegraphics[height=0.5cm]{diri07.eps}), \pi_{0,\diamond}(\includegraphics[height=0.7cm]{diri09.eps},\includegraphics[height=0.5cm]{diri07.eps})$ are defined  as in Appendix B. To simplify the arguments below, we assume that $\mathcal{F}{W}(0)=0$ and restrict ourselves to the flow of $\int_{\mathbb{T}^3} u(x)dx=0$. Furthermore, we use $\includegraphics[height=0.7cm]{diri05.eps}$ and $\includegraphics[height=0.7cm]{diri10.eps}$ to denote $-\rho_\varepsilon*\Phi_2$ and $\rho_\varepsilon*K$, respectively. We summarise the graph notations after the introduction.   We also introduce the following:
 $$\pi_{0,\diamond}(\includegraphics[height=0.7cm]{diri05.eps},\includegraphics[height=0.5cm]{diri07.eps})
 :=\pi_{0}(\includegraphics[height=0.7cm]{diri05.eps},\includegraphics[height=0.5cm]{diri07.eps})-3(\tilde{C}_{1}^\varepsilon+\tilde{\varphi}^\varepsilon) \includegraphics[height=0.5cm]{diri01.eps},$$
 $$\pi_{0,\diamond}(\includegraphics[height=0.7cm]{diri10.eps},\includegraphics[height=0.5cm]{diri07.eps})
 :=\pi_{0}(\includegraphics[height=0.7cm]{diri10.eps},\includegraphics[height=0.5cm]{diri07.eps})-(\tilde{C}_{1}^\varepsilon+\tilde{\varphi}^\varepsilon),$$
 with $$\tilde{C}_{1}^\varepsilon=2^{-7}\int\int\frac{g(\varepsilon k_1)g(\varepsilon k_2)g(\varepsilon k_{[12]})}{|k_1|^2|k_2|^2(|k_1|^2+|k_2|^2+|k_{[12]}|^2)\pi^6}dk_1dk_2,$$
 and
 $$\tilde{\varphi}^\varepsilon(t)=-2^{-7}\int\int\frac{e^{-t\pi^2(|k_1|^2+|k_2|^2+|k_{[12]}|^2)}g(\varepsilon k_1)g(\varepsilon k_2)g(\varepsilon k_{[12]})}{|k_1|^2|k_2|^2(|k_1|^2+|k_2|^2+|k_{[12]}|^2)\pi^6}dk_1dk_2.$$
Here $k_{[12]}=k_1+k_2$ and the integral is on the set $\mathbb{Z}^3\backslash \{0\}$.

We also define
 $$\aligned\delta \bar{C}_W^{\varepsilon}(T):=&\sup_{t\in[0,T]}\big[\|\pi_0(\includegraphics[height=0.7cm]{diri05.eps},\includegraphics[height=0.5cm]{diri02.eps})
 -\pi_0(\includegraphics[height=0.7cm]{diri04.eps},\includegraphics[height=0.5cm]{diri02.eps})\|_{-2\kappa}
 +\|\pi_{0,\diamond}(\includegraphics[height=0.7cm]{diri05.eps},\includegraphics[height=0.5cm]{diri07.eps})
 -\pi_{0,\diamond}(\includegraphics[height=0.7cm]{diri04.eps},\includegraphics[height=0.5cm]{diri07.eps})\|_{-\frac{1}{2}-2\kappa}
 \\&+\|\pi_{0,\diamond}(\includegraphics[height=0.7cm]{diri10.eps},\includegraphics[height=0.5cm]{diri07.eps})
 -\pi_{0,\diamond}(\includegraphics[height=0.7cm]{diri09.eps},\includegraphics[height=0.5cm]{diri07.eps})\|_{-2\kappa}
 +\|\includegraphics[height=0.5cm]{diri01.eps}-\includegraphics[height=0.5cm]{diri02.eps}\|_{-\frac{1}{2}-2\kappa}\\&
 +\|\includegraphics[height=0.5cm]{diri06.eps}-\includegraphics[height=0.5cm]{diri07.eps}\|_{-1-2\kappa}
 +\|\includegraphics[height=0.7cm]{diri03.eps}-\includegraphics[height=0.7cm]{diri04.eps}\|_{\frac{1}{2}-2\kappa}+\|\pi_{0}( \includegraphics[height=0.7cm]{diri03.eps} ,\includegraphics[height=0.5cm]{diri01.eps} )-\pi_0(\includegraphics[height=0.7cm]{diri04.eps},\includegraphics[height=0.5cm]{diri02.eps})\|_{-2\kappa}
 \\&+\|\pi_{0,\diamond}( \includegraphics[height=0.7cm]{diri03.eps} ,\includegraphics[height=0.5cm]{diri06.eps})
 -\pi_{0,\diamond}(\includegraphics[height=0.7cm]{diri04.eps},\includegraphics[height=0.5cm]{diri07.eps})\|_{-\frac{1}{2}-2\kappa}+\|\pi_{0,\diamond}
(\includegraphics[height=0.7cm]{diri08.eps} ,\includegraphics[height=0.5cm]{diri06.eps})-\pi_{0,\diamond}(\includegraphics[height=0.7cm]{diri09.eps},\includegraphics[height=0.5cm]{diri07.eps})
\|_{-2\kappa}\big]\\&+\|\includegraphics[height=0.7cm]{diri03.eps}-\includegraphics[height=0.7cm]{diri04.eps}\|_{C^{\frac{1}{8}}_T
\mathcal{C}^{\frac{1}{4}-2\kappa}}.\endaligned$$
By Appendix B we can find a subsequence of $\varepsilon$ going to zero such that for any $T>0$ $\lim_{\varepsilon\rightarrow0}\delta \bar{C}_W^\varepsilon(T)=0, \lim_{\varepsilon\rightarrow0}\int_0^T\includegraphics[height=0.5cm]{diri11.eps}ds$ exists $P$-a.s.. Here and in the following for simplicity we still use the notation $\varepsilon$ to denote this subsequence.
Set $$\Omega_0=\{ \lim_{\varepsilon\rightarrow0}\delta \bar{C}_W^\varepsilon(T)=0, C_W(T)<\infty, \lim_{\varepsilon\rightarrow0}\int_0^T\includegraphics[height=0.5cm]{diri11.eps}ds \textrm{ exists, for any } T>0\}.$$
Then $P(\Omega_0)=1$.
\vskip.10in
\th{Lemma 5.1} $\Phi$ satisfies (5.1) in the analytically weak sense on $\Omega_0$.

\proof First we prove the following:
$$\aligned \lim_{\varepsilon\rightarrow0}\int_0^t[(\rho_\varepsilon\ast\Phi)^3-(3\bar{C}_0^\varepsilon\rho_\varepsilon*\Phi-9\tilde{C}_{1}^\varepsilon\Phi-m\Phi) ]ds
=\lim_{\varepsilon\rightarrow0}\int_0^t[(\bar{\Phi}^\varepsilon)^3-(3\bar{C}_0^\varepsilon-9\bar{C}_{1}^\varepsilon-m)\bar{\Phi}^\varepsilon ]ds\endaligned.\eqno(5.3)$$

In fact,
$$\aligned&\int_0^t[(\bar{\Phi}^\varepsilon)^3-(3\bar{C}_0^\varepsilon-9\bar{C}_{1}^\varepsilon-m) \bar{\Phi}^\varepsilon ]ds\\=&\int_0^t[(\bar{\Phi}^\varepsilon_3)^3 +3(\includegraphics[height=0.5cm]{diri02.eps}-\includegraphics[height=0.7cm]{diri04.eps}) (\bar{\Phi}^\varepsilon_3 )^2+(3(\includegraphics[height=0.7cm]{diri04.eps})^2-6\includegraphics[height=0.5cm]{diri02.eps}  \includegraphics[height=0.7cm]{diri04.eps}) \bar{\Phi}^\varepsilon_3+3\includegraphics[height=0.5cm]{diri07.eps}\diamond \bar{\Phi}^\varepsilon_3\\&+\includegraphics[height=0.5cm]{diri11.eps}+3\includegraphics[height=0.5cm]{diri02.eps} (\includegraphics[height=0.7cm]{diri04.eps} )^{2}-(\includegraphics[height=0.7cm]{diri04.eps})^3-3\includegraphics[height=0.5cm]{diri07.eps}\diamond\includegraphics[height=0.7cm]{diri04.eps} -(9\bar{\varphi}^\varepsilon-m) \bar{\Phi}^\varepsilon]ds,\endaligned\eqno(5.4)$$
and
$$\aligned&\int_0^t[(\rho_\varepsilon*\Phi)^3-(3\bar{C}_0^\varepsilon\rho_\varepsilon*\Phi-9\tilde{C}_{1}^\varepsilon\Phi-m\Phi) ]ds\\=&\int_0^t[(\rho_\varepsilon*\Phi_3)^3 +3(\includegraphics[height=0.5cm]{diri02.eps}-\includegraphics[height=0.7cm]{diri05.eps}) (\rho_\varepsilon*\Phi_3 )^2+(3(\includegraphics[height=0.7cm]{diri05.eps})^2-6\includegraphics[height=0.5cm]{diri02.eps}  \includegraphics[height=0.7cm]{diri05.eps}) \rho_\varepsilon*\Phi_3+3\includegraphics[height=0.5cm]{diri07.eps}\diamond \rho_\varepsilon*\Phi_3\\&+\includegraphics[height=0.5cm]{diri11.eps}+3\includegraphics[height=0.5cm]{diri02.eps} (\includegraphics[height=0.7cm]{diri05.eps} )^{2}-(\includegraphics[height=0.7cm]{diri05.eps})^3-3\includegraphics[height=0.5cm]{diri07.eps}\diamond\includegraphics[height=0.7cm]{diri05.eps} -(9\tilde{\varphi}^\varepsilon-m) {\Phi}]ds,\endaligned\eqno(5.5)$$
where
$$\aligned\includegraphics[height=0.5cm]{diri07.eps} \diamond\includegraphics[height=0.7cm]{diri05.eps}:=&\includegraphics[height=0.5cm]{diri07.eps} \includegraphics[height=0.7cm]{diri05.eps}-3({\tilde{C}}_{1}^\varepsilon+{\tilde{\varphi}}^\varepsilon) \includegraphics[height=0.5cm]{diri01.eps},\\ \includegraphics[height=0.5cm]{diri07.eps}\diamond \rho_\varepsilon*\Phi_3:=&\rho_\varepsilon*\Phi_3\includegraphics[height=0.5cm]{diri07.eps}+3({\tilde{C}}_{1}^\varepsilon+{\tilde{\varphi}}^\varepsilon) (-\includegraphics[height=0.7cm]{diri03.eps}+{\Phi}_3),\endaligned$$
and the other terms containing $\diamond$ and  $\bar{\varphi}^\varepsilon$ are defined  in Appendix B and $\Phi_3$ satisfies equation (2.1).  Now we only need to prove that each term converges. First we check the relations between  $\includegraphics[height=0.7cm]{diri05.eps}$, $\rho_\varepsilon*\Phi_3$ and $\includegraphics[height=0.7cm]{diri04.eps}, \bar{\Phi}_3^\varepsilon$.  We have that on $\Omega_0$  for any $T>0$
and  $\epsilon>0$ small enough
$$\sup_{t\in [0,T]}\|\includegraphics[height=0.7cm]{diri05.eps}-\includegraphics[height=0.7cm]{diri04.eps}\|_{\frac{1}{2}-2\kappa-\epsilon}
\leq \sup_{t\in [0,T]}(\|\includegraphics[height=0.7cm]{diri03.eps}-\includegraphics[height=0.7cm]{diri04.eps}\|_{\frac{1}{2}-2\kappa}
+\|\includegraphics[height=0.7cm]{diri05.eps}-\includegraphics[height=0.7cm]{diri03.eps}\|_{\frac{1}{2}-2\kappa-\epsilon})\rightarrow0.$$
Now we consider $\rho_\varepsilon*\Phi_3-\bar{\Phi}_3$. We define $\bar{C}_W^\varepsilon(T, \omega)$ for (5.2) similarly as $C_W(T,\omega)$ in (2.2) and
 we have that for $\omega\in\Omega_0$, there exists a constant $C_1(T,\omega)$ such that $\bar{C}_W^\varepsilon(T, \omega)\leq C_1(T,\omega)$ for the subsequence of $\varepsilon$. Since $\bar{\Phi}_3^\varepsilon$ satisfies a similar equation as $\Phi_3$, by a similar argument as in Proposition 2.1 we obtain that $$\sup_{t\in[0,T]}[t^{\frac{\gamma+z+\kappa}{2}}(\|\Phi_3\|_\gamma+\|\bar{\Phi}^\varepsilon_3\|_\gamma)+t^{\frac{\frac{1}{2}+z+5\kappa}{2}}
 (\|\Phi_3\|_{\frac{1}{2}+4\kappa}+\|\bar{\Phi}^\varepsilon_3\|_{\frac{1}{2}+4\kappa})]\leq C(T,\omega,\|\Phi_0\|_{-z}).$$Then  a similar argument as in Proposition 3.1 yields that on $\Omega_0$
 $$\sup_{t\in[0,T]}[t^{\frac{\gamma+z+\kappa}{2}}\|\Phi_3-\bar{\Phi}^\varepsilon_3\|_{\gamma}+t^{\frac{\frac{1}{2}+z+5\kappa}{2}}
 \|\Phi_3-\bar{\Phi}^\varepsilon_3\|_{\frac{1}{2}+4\kappa}+t^{\frac{3(\gamma+z+\kappa)}{2}}\|\Phi^\sharp-\bar{\Phi}^{\varepsilon,\sharp}\|_{1+3\kappa}]\rightarrow0,$$ which combined with the fact the $\|\rho_\varepsilon*\Phi_3-\Phi_3\|_{\beta-\kappa}\lesssim \varepsilon^{\frac{\kappa}{2}}\|\Phi_3\|_\beta$ implies that on $\Omega_0$ for $\epsilon>0$ small enough $$\sup_{t\in[0,T]}[t^{\frac{\gamma+z+\kappa}{2}}\|\rho_\varepsilon*\Phi_3-\bar{\Phi}^\varepsilon_3\|_{\gamma-\epsilon}+t^{\frac{\frac{1}{2}+z+5\kappa}{2}}
 \|\rho_\varepsilon*\Phi_3-\bar{\Phi}^\varepsilon_3\|_{\frac{1}{2}+4\kappa-\epsilon}
 +t^{\frac{3(\gamma+z+\kappa)}{2}}\|\rho_\varepsilon*\Phi^\sharp-\bar{\Phi}^{\varepsilon,\sharp}\|_{1+3\kappa-\epsilon}]\rightarrow0.\eqno(5.6)$$ Hence by Lemma A.2 we obtain that the terms which do not need to be renormalized in (5.4) and (5.5) converge. Now we concentrate on the renormalization terms. For the renormalized terms $\includegraphics[height=0.5cm]{diri02.eps}\includegraphics[height=0.7cm]{diri04.eps}, \includegraphics[height=0.5cm]{diri07.eps}\diamond\includegraphics[height=0.7cm]{diri04.eps}, \includegraphics[height=0.5cm]{diri02.eps}(\includegraphics[height=0.7cm]{diri04.eps})^2$ by Lemmas A.2, A.3 it is sufficient to consider the following terms: Since $\delta\bar{C}_W^\varepsilon\rightarrow0$ on $\Omega_0$, we have on $\Omega_0$
$$\pi_0(\includegraphics[height=0.7cm]{diri05.eps},\includegraphics[height=0.5cm]{diri02.eps})
-\pi_0(\includegraphics[height=0.7cm]{diri04.eps},\includegraphics[height=0.5cm]{diri02.eps})\rightarrow0\quad \textrm{in } C_T\mathcal{C}^{-2\kappa},$$
and$$\pi_{0,\diamond}(\includegraphics[height=0.7cm]{diri05.eps},\includegraphics[height=0.5cm]{diri07.eps})
-\pi_{0,\diamond}(\includegraphics[height=0.7cm]{diri04.eps},\includegraphics[height=0.5cm]{diri07.eps})\rightarrow0\quad \textrm{in } C_T\mathcal{C}^{-\frac{1}{2}-2\kappa}.$$
Now we focus on the convergence of $\includegraphics[height=0.5cm]{diri07.eps}\diamond\rho_\varepsilon*\Phi_3 $. It is sufficient to consider $\pi_{0,\diamond}(\rho_\varepsilon*\Phi_3,\includegraphics[height=0.5cm]{diri07.eps})
:=\pi_{0}(\rho_\varepsilon*\Phi_3,\includegraphics[height=0.5cm]{diri07.eps})+3(-\includegraphics[height=0.7cm]{diri03.eps}
+\Phi_3)\pi_{0,\diamond}(\includegraphics[height=0.7cm]{diri10.eps},\includegraphics[height=0.5cm]{diri07.eps})$. We have $$\Phi_3=-3\pi_<(-\includegraphics[height=0.7cm]{diri03.eps}+\Phi_3,\includegraphics[height=0.7cm]{diri08.eps})+\Phi^\sharp.$$ Then we obtain that
$$\pi_0(\rho_\varepsilon*\Phi_3,\includegraphics[height=0.5cm]{diri07.eps})=-3\pi_0(\rho_\varepsilon*\pi_<(-\includegraphics[height=0.7cm]{diri03.eps}
+\Phi_3,\includegraphics[height=0.7cm]{diri08.eps}),\includegraphics[height=0.5cm]{diri07.eps})
+\pi_0(\rho_\varepsilon*\Phi^\sharp,\includegraphics[height=0.5cm]{diri07.eps}).$$
For the second term we can easily obtain the convergence by (5.6). For the first term we have
$$\aligned &\pi_0(\rho_\varepsilon*\pi_<(-\includegraphics[height=0.7cm]{diri03.eps}+\Phi_3,\includegraphics[height=0.7cm]{diri08.eps}),
\includegraphics[height=0.5cm]{diri07.eps})
\\=&\pi_0(\rho_\varepsilon*\pi_<(-\includegraphics[height=0.7cm]{diri03.eps}+\Phi_3,\includegraphics[height=0.7cm]{diri08.eps})
,\includegraphics[height=0.5cm]{diri07.eps})-\pi_0(\pi_<(-\includegraphics[height=0.7cm]{diri03.eps}+\Phi_3,\includegraphics[height=0.7cm]{diri10.eps}),
\includegraphics[height=0.5cm]{diri07.eps})
\\&+C(-\includegraphics[height=0.7cm]{diri03.eps}+\Phi_3,\includegraphics[height=0.7cm]{diri10.eps},\includegraphics[height=0.5cm]{diri07.eps})
+(-\includegraphics[height=0.7cm]{diri03.eps}+\Phi_3)\pi_0(\includegraphics[height=0.7cm]{diri10.eps},\includegraphics[height=0.5cm]{diri07.eps}),
\endaligned$$
where the first two terms converge to zero as $\varepsilon\rightarrow0$ by  Lemma A.5 and the third  term converges to the corresponding term by Lemma A.3 and the last term should be renormalized and converges to the corresponding term on $\Omega_0$. Since $\lim_{\varepsilon\rightarrow0}\int_0^t\includegraphics[height=0.5cm]{diri11.eps}ds$ exists, combining the above arguments (5.3) follows. Moreover, on $\Omega_0$ we know that for any $t>0$,

$$\bar{\Phi}^\varepsilon(t)= \Phi_0+\int_0^t\Delta \bar{\Phi}^\varepsilon ds-\int_0^t[(\bar{\Phi}^\varepsilon)^3-(3 \bar{C}_0^\varepsilon-9\bar{C}_1^\varepsilon-m)\bar{\Phi}^\varepsilon ]ds +\rho_\varepsilon*W(t).$$ Then taking the limit on both sides we obtain the result. $\hfill\Box$

\vskip.10in

\no \emph{Proof of Theorem 1.2} The idea is to prove that the drift term in (5.1) is the zero-energy part in the Fukushima decomposition (cf. [FOT94, Theorem 5.2.2]).
In the proof we take the space of continuous paths $C([0,\infty);E)$ as the sample paths $\bar{\Omega}$ and we denote the $t$-th coordinate of the path $\omega$ by $\bar{X}_t(\omega)$.
For $t\in[0,\infty)$ let $(\bar{\mathcal{F}}_t)$ be the natural filtration for $\bar{X}$ given in [MR92, Chapter IV, (1.7)]. Set $\bar{\mathcal{F}}:=\cup_{t\geq0}\bar{\mathcal{F}}_t$ and define on $\bar{\Omega}$ $$P^x(\bar{X}\in A):=P(\Phi(x)\in A),$$
for $A\in\mathcal{B}(\bar{\Omega})$. Here $\Phi$ on the right hand side is the solution from Section 2 starting from $x$. Under $P^x$, $\bar{X}$ is the solution to (1.2) starting from $x$. Let $\theta$ be the associated shift operator. By  Theorem D.4 the  (Markov) diffusion process
$(\Omega, \bar{\mathcal{F}}, (\bar{\mathcal{F}}_t)_{t>0}, \theta_t, \bar{X}, P^x )_{x\in E}$ is properly associated with $(\mathcal{E}, D(\mathcal{E}))$.
Define
$$\Omega_1:=\{\omega:\lim_{\varepsilon\rightarrow0}\int_0^{\cdot}\langle(\rho_\varepsilon\ast\bar{X})^3-(3\bar{C}_0^\varepsilon\rho_\varepsilon*\bar{X}
-9\tilde{C}_{1}^\varepsilon\bar{X}-m\bar{X}),\varphi\rangle dr \textrm{ exists in } C([0,\infty);\mathbb{R}) , \forall\varphi\in\mathcal{D}\}.$$
and for $\varphi\in C^\infty(\mathbb{T}^3)$, $$H_t^\varphi:=\left\{\begin{array}{ll}\lim_{\varepsilon\rightarrow0}\int_0^t\langle(\rho_\varepsilon\ast\bar{X})^3
-(3\bar{C}_0^\varepsilon\rho_\varepsilon*\bar{X}-9\tilde{C}_{1}^\varepsilon\bar{X}-m\bar{X}), \varphi\rangle ds,&\ \ \ \ \textrm{ for } \omega\in\Omega_1
\\0&\ \ \ \ \textrm{ otherwise } .\end{array}\right.
$$
Now we would like to check that $H_t^\varphi$ is an additive functional (AF) in the sense of [FOT94, Section 5.1]:

i) It's obvious that  $H_t^\varphi$ is $\bar{\mathcal{F}}_t$-measurable;

ii) For $\omega\in\Omega$,  $H_\cdot^\varphi(\omega)$ is continuous, $H_0(\omega)=0$. Since $P^x(\bar{X}\in C([0,\infty);\mathcal{C}^{-z}))=1$ for $x\in \mathcal{C}^{-z}$ and $\mu(\mathcal{C}^{-z})=1$, it is sufficient to check that for $x\in \mathcal{C}^{-z}$ $P^x(\Omega_1)=1$, $\theta_t\Omega_1\subset \Omega_1$, and for $\omega\in\Omega_1$
$$H_{t+s}^\varphi(\omega)=H_{t}^\varphi(\omega)+H_{s}^\varphi(\theta_t\omega).\eqno(5.7)$$

$P(\Omega_0)=1$ implies that $P^x(\Omega_1)=1$ by Lemma 5.1. Since $\bar{X}(t+s)=\bar{X}(s)\circ \theta_t$, we can easily deduce that $\theta_t\Omega_1\subset \Omega_1$ and that
$$\aligned&\int_0^{t+s}\langle(\rho_\varepsilon\ast\bar{X})^3-(3\bar{C}_0^\varepsilon\rho_\varepsilon*\bar{X}-9\tilde{C}_{1}^\varepsilon\bar{X}-m\bar{X}), \varphi\rangle dr\\=&\int_0^{t}\langle(\rho_\varepsilon\ast\bar{X})^3-(3\bar{C}_0^\varepsilon\rho_\varepsilon*\bar{X}-9\tilde{C}_{1}^\varepsilon\bar{X}-m\bar{X}), \varphi\rangle dr\\+&\int_0^{s}\langle(\rho_\varepsilon\ast\bar{X})^3-(3\bar{C}_0^\varepsilon\rho_\varepsilon*\bar{X}-9\tilde{C}_{1}^\varepsilon\bar{X}-m\bar{X}), \varphi\rangle dr\circ\theta_t,\endaligned$$
which implies that (5.7) holds for $\omega\in\Omega_1$.

Now we know that $H_t^\varphi$ is an AF. Define
$$M^\varphi_t:=\langle\bar{X}(t)- \bar{X}(0),\varphi\rangle-\int_0^t\langle \bar{X},\Delta \varphi\rangle ds+H_t^\varphi.$$
We know that $M^\varphi$ is also an AF. Moreover, by Lemma 5.1 we have
$$E^x M^\varphi_t=0,\quad E^x (M^\varphi_t)^2=|\varphi|^2t<\infty,$$
which implies that $M^\varphi$ is also a martingale additive functional (MAF) in the sense of [FOT94, Chapter V]. Here $|\cdot|$ denotes the $L^2$-norm.

 Let us fix an arbitrary $T>0$ and consider the  space $\Omega_T$ of all continuous paths from $[0,T]$ to $E$.
 We introduce the time reversal operator $r_T$ on $\Omega_T$ defined by
 $$r_T\omega(t)=\omega(T-t),\quad 0\leq t\leq T, \omega\in \Omega_T.$$
 By [FOT94, Lemma 5.7.1] and the symmetry of the semigroup $\bar{P}_t$ we have that for any $\bar{\mathcal{F}}_T$-measurable set $A$ on $\Omega_T$
$$P^\mu(r_T\omega\in A)=P^\mu(A),\eqno(5.8)$$
where $P^\mu=\int P^x\mu(dx)$.
Now we have
$$\langle\bar{X}(t)- \bar{X}(0),\varphi\rangle=M^\varphi_t+\bar{H}_t^\varphi\quad P^\mu-a.s.,$$
with $\bar{H}_t^\varphi=\int_0^t\langle \bar{X},\Delta \varphi\rangle ds-H_t^\varphi$.
By (5.8) we have  for $0\leq t\leq T$
$$\langle\bar{X}(T-t)- \bar{X}(T),\varphi\rangle=M^\varphi_t(r_T)+\bar{H}_t^\varphi(r_T)\quad P^\mu-a.s..\eqno(5.9)$$
Moreover, under $P^\mu$,
$$\aligned\bar{H}_t^\varphi(r_T)=&\int_0^t\langle \bar{X}\circ r_T,\Delta \varphi\rangle ds-\lim_{\varepsilon\rightarrow0}\int_0^t\langle((\rho_\varepsilon\ast\bar{X})\circ r_T)^3-(3\bar{C}_0^\varepsilon\rho_\varepsilon*\bar{X}-9\tilde{C}_{1}^\varepsilon\bar{X}-m\bar{X})\circ r_T, \varphi\rangle ds\\=&\int_{T-t}^T\langle \bar{X},\Delta \varphi\rangle ds-\lim_{\varepsilon\rightarrow0}\int_{T-t}^T\langle(\rho_\varepsilon\ast\bar{X})^3-(3\bar{C}_0^\varepsilon\rho_\varepsilon*\bar{X}
-9\tilde{C}_{1}^\varepsilon\bar{X}-m\bar{X}), \varphi\rangle ds\\=&\bar{H}_T^\varphi-\bar{H}_{T-t}^\varphi.\endaligned\eqno(5.10)$$
By (5.9), (5.10) we have
$$M^\varphi_t(r_T)=\langle\bar{X}(T-t)- \bar{X}(T),\varphi\rangle-\bar{H}_T^\varphi+\bar{H}_{T-t}^\varphi,$$
which implies that
$$\aligned &M^\varphi_{T-t}(r_T)-M^\varphi_{T}(r_T)
\\=&\langle\bar{X}(t)- \bar{X}(T),\varphi\rangle-\bar{H}_T^\varphi+\bar{H}_{t}^\varphi-\langle\bar{X}(0)- \bar{X}(T),\varphi\rangle+\bar{H}_T^\varphi
\\=&2\langle\bar{X}(t)- \bar{X}(0),\varphi\rangle-M^\varphi_{t}.\endaligned$$
Now we know that
$$\langle\bar{X}(t)- \bar{X}(0),\varphi\rangle=\frac{1}{2}(M^\varphi_{t}-M^\varphi_{t}\circ r_t)\quad P^\mu-a.s. \forall t>0.$$
By [F95, Theorem 2.2] we have that $M^\varphi\equiv M^{[\varphi]}$, where $M^{[\varphi]}$ is the MAF from the Fukushima decomposition for $\langle \cdot,\varphi\rangle$ (see [FOT94, Section 5.2]. Hence, we have that $\bar{H}_t^\varphi=N_t^{[\varphi]}$ is the associated zero-energy additive functional (NAF), which implies that $\Phi$ is a Dirichlet process.

Now for  $f=f_1(\langle \cdot, l_1\rangle,\langle \cdot, l_2\rangle,...,\langle \cdot, l_k\rangle)$ with $l_i, f_1$ smooth, denote the MAF in the Fukushima decomposition associated with $\langle\cdot,l_i\rangle$ by $M^{l_i}$.
By It\^{o}'s formula for Dirichlet process in [CFKZ08, Theorem 4.7] and [N85, Theorem 4.1], we have
$$\aligned f(\bar{X}(t))-f(\bar{X}(0))=&\sum_{i=1}^k \int_0^t \partial_if(\bar{X}(s))d M^{l_i}_{t}+ \sum_{i=1}^k \int_0^t \partial_if(\bar{X}(s))d \bar{H}^{l_i}_{t} +\frac{1}{2}\sum_{i,j=1}^k\int_0^t \partial_{ij}f(\bar{X}(s))\langle l_i,l_j\rangle ds\\:=&\sum_{i=1}^k \int_0^t \partial_if(\bar{X}(s))d M^{l_i}_{t}+\bar{H}_t^f,\endaligned$$
where $\partial_if:=\partial_if_1(\langle \cdot, l_1\rangle,\langle \cdot, l_2\rangle,...,\langle \cdot, l_k\rangle)$ and the stochastic integral $\int_0^t \partial_if(\bar{X}(s))d \bar{H}^{l_i}_{t} $ w.r.t. NAF is defined in [CFKZ08].
We know that $\sum_{i=1}^k \int_0^t \partial_if(\bar{X}(s))d M^{l_i}_{t}$ is an MAF and $\bar{H}_t^f$ is an NAF, which implies that
$$\sum_{i=1}^k \int_0^t \partial_if(\bar{X}(s))d M^{l_i}_{t}\equiv M_t^{[f]},\eqno(5.11)$$
where $M_t^{[f]}$ is the MAF obtained in the Fukushima decomposition.

By (5.11) we know that $$\mathcal{E}(f,f)=e(M_t^{[f]}):=\lim_{t\downarrow0}\frac{1}{2t}E^\mu(M_t^{[f]})^2=\frac{1}{2}\int|Df|^2d\mu.$$
Then for $g\in\mathcal{F}C_b^\infty$ we can use the above $f$'s to approximate it and obtain $\mathcal{E}(g,g)=\frac{1}{2}\int |Dg|^2d\mu$.$\hfill\Box$

\vskip.10in
\th{Remark 5.2} From the above proof we can  check that $\Phi$ starting from $\mu$ is an energy solution in the sense that $(\Phi,N)_{0\leq t\leq T}$ has continuous paths in $E$ such that

i) the law of $\Phi$ is  $\mu$ for all $t\in[0, T]$;

ii) for any test function $\varphi\in C^\infty(\mathbb{T}^3)$ the process $t\rightarrow N_t$ is a.s. of zero quadratic variation, $N_0(\varphi)=0$
and the pair $(\Phi(\varphi),N(\varphi))_{0\leq t\leq T}$ satisfies the equation
$$\langle\Phi_t,\varphi\rangle=\langle\Phi_0,\varphi\rangle+\int_0^t\langle\Phi_s,\Delta\varphi\rangle ds+\langle N_t,\varphi\rangle+\langle M_t,\varphi\rangle,$$
where $(\langle M_t,\varphi\rangle)_{0\leq t\leq T}$ is a martingale with respect to the filtration generated by $(\Phi,N)_{0\leq t\leq T}$ with
quadratic variation $|\varphi|^2t$.

iii) the reversed processes $\hat{\Phi}_t = \Phi_{T-t}, \hat{N}_t = N_T-N_{T-t}$ satisfies the same equation with the associated martingale $\hat{M}_t$ with respect
to its own filtration and the quadratic variation of $\hat{M}$ is also $|\varphi|^2t$.

iv) $ N_t=-\lim_{\varepsilon\rightarrow0}\int_0^t[(\rho_\varepsilon\ast\Phi)^3-(3\bar{C}_0^\varepsilon\rho_\varepsilon*\Phi-9\tilde{C}_{1}^\varepsilon\Phi-m\Phi) ]ds$ a.s. with $\bar{C}_0^\varepsilon, \tilde{C}_{1}^\varepsilon$ introduced at the beginning of Section 5.

\vskip.10in
\no\emph{Proof of Theorem 1.3} By Theorem 1.2 we know that $(\bar{\mathcal{E}},\mathcal{F}C_b^\infty)$ is a well-defined symmetric bilinear form. Since the Dirichlet form $(\mathcal{E}, D(\mathcal{E}))$ is an extension of $(\bar{\mathcal{E}}, \mathcal{F}C_b^\infty)$, it is obvious that  $(\bar{\mathcal{E}}, \mathcal{F}C_b^\infty)$ is closable. We denote its closure by $(\bar{\mathcal{E}}, D(\bar{\mathcal{E}}))$. Then by similar arguments
as in [MR92, Chapter II Proposition 3.5] we obtain that for $u\in D(\bar{\mathcal{E}})$, $v=u\vee0\wedge1\in D(\bar{\mathcal{E}})$ and $\bar{\mathcal{E}}(v,v)\leq \bar{\mathcal{E}}(u,u)$. Moreover, by similar arguments as in the proof of [MR92, Chapter IV Proposition 4.2] i) in Definition D.1 follows, which implies
 $(\bar{\mathcal{E}}, D(\bar{\mathcal{E}}))$ is a quasi-regular Dirichlet form (cf. Definition D.1). Then existence of the Markov process follows from Theorem D.4.
$\hfill\Box$

\vskip.10in
\no\emph{Proof of Corollary 1.5} By general theory of Markov semigroup and Dirichlet form (cf. [W05]) we know the following Poincar\'{e} inequality holds:
$$\mu(f^2)\leq C\mathcal{E}(f,f)+\mu(f)^2, \quad f\in D(\mathcal{E})\eqno(5.12)$$
for some $C>0$. In the following we follows essentially the same argument from [W05, Section 1.2] to deduce the last result. Since $$\|x\|_E^2=\sum_k\lambda_k\langle x,\hat{e}_k\rangle^2,$$
where $\lambda_k\in\mathbb{R}$ satisfies $\lambda_k\rightarrow0, k\rightarrow\infty$ and $\{\hat{e}_k\}$ is a real smooth eigenbasis on $L^2(\mathbb{T}^3)$. We first prove that for $r\geq0, n\in\mathbb{N}$, $f_n(\cdot):=e^{\frac{r}{2}((\sum_{k}\lambda_k\langle \cdot,\hat{e}_k\rangle^2+1)^{\frac{1}{2}}\wedge n)}\in D(\mathcal{E})$.
By approximation we can easily check that
$f_{n,N}:=e^{\frac{r}{2}((\sum_{|k|_\infty\leq N}\lambda_k\langle \cdot,\hat{e}_k\rangle^2+1)^{\frac{1}{2}}\wedge n)}\in D(\mathcal{E})$.
Moreover, by direct computation we know that
$$\mathcal{E}_1(f_{n,N}-f_n,f_{n,N}-f_n)\rightarrow0,\quad N\rightarrow\infty,$$
with $\mathcal{E}_1(\cdot,\cdot):=\mathcal{E}(\cdot,\cdot)+(\cdot,\cdot)_{L^2(E;\mu)}$. 
We also have
$$\mathcal{E}(f_{n,N},f_{n,N})\leq \frac{r^2}{4}\int f_{n,N}^2d\mu, $$
which implies the following by letting $N\rightarrow\infty$
 $$\mathcal{E}(f_n,f_n)\leq \frac{r^2}{4}\int f_n^2d\mu. $$
 Let $h_n(r):=\mu(f_n^2)$. 
  By (5.12) we know that 
 $$h_n(r)\leq \frac{Cr^2}{4}h_n(r)+h_n(r/2)^2.$$
Thus, for any $r\in (0,2/\sqrt{C})$ we have
$$h_n(r)\leq \frac{4}{4-Cr^2}h_n(r/2)^2.\eqno(5.13)$$
Next, for any $m>0$, let $p_m=\mu(x:(\sum_{k}\lambda_k\langle x,\hat{e}_k\rangle^2+1)^{1/2}\geq m)$. We have 
$$h_n(r/2)^2\leq \bigg[e^{mr/2}+\mu(1_{\{(\sum_{k}\lambda_k\langle x,\hat{e}_k\rangle^2+1)^{1/2}\geq m\}}f_n)\bigg]^2\leq 2e^{mr}+2p_mh_n(r).$$
Substituting this into (5.13) we have 
$$h_n(r)\leq \frac{8}{4-Cr^2}e^{mr}+\frac{8}{4-Cr^2}p_mh_n(r), \quad 0<r<2/\sqrt{C}.$$
By Lemma 2.3 we know that $p_m\rightarrow0$ as $m\rightarrow\infty$, which implies that there exists $m_0>0$ such that 
$\frac{8p_{m_0}}{4-Cr^2}\leq \frac{1}{2}$. Therefore, 
$$h_n(r)\leq \frac{16}{4-Cr^2}e^{m_0r}.$$
Letting $n\rightarrow\infty$ we arrive at 
$$\int e^{r\|x\|_E}\mu(dx)\leq \int e^{r(\sum_{k}\lambda_k\langle x,\hat{e}_k\rangle^2+1)^{\frac{1}{2}}}\mu(dx)<\infty, \quad r\in (0,2/\sqrt{C}).$$$\hfill\Box$

\vskip.20in

\no\textbf{Acknowledgments}
\vskip.10in
We are very grateful to Professor Zhenqing Chen for pointing out reference [F95] to us and helpful discussions. We would also like to thank Professor Michael R\"{o}ckner and Professor Lorenzo Zambotti for their encouragement and suggestions for this work.
\vskip.20in
\no\textbf{Appendix A: Besov spaces and paraproduct}
\vskip.10in
In this appendix we recall the definitions and some properties of Besov spaces and paraproducts. For a general introduction to these theories we refer to [BCD11, GIP15].
First we introduce the following notations. The space of real valued infinitely differentiable functions of compact support is denoted by $\mathcal{D}(\mathbb{R}^d)$ or $\mathcal{D}$. The space of Schwartz functions is denoted by $\mathcal{S}(\mathbb{R}^d)$. Its dual, the space of tempered distributions is denoted by $\mathcal{S}'(\mathbb{R}^d)$.

 Let $\chi,\theta\in \mathcal{D}$ be nonnegative radial functions on $\mathbb{R}^d$, such that

i. the support of $\chi$ is contained in a ball and the support of $\theta$ is contained in an annulus;

ii. $\chi(z)+\sum_{j\geq0}\theta(2^{-j}z)=1$ for all $z\in \mathbb{R}^d$.

iii. $\textrm{supp}(\chi)\cap \textrm{supp}(\theta(2^{-j}\cdot))=\emptyset$ for $j\geq1$ and $\textrm{supp}(\theta(2^{-i}\cdot))\cap \textrm{supp}(\theta(2^{-j}\cdot))=\emptyset$ for $|i-j|>1$.

We call such a pair $(\chi,\theta)$ a dyadic partition of unity, and for the existence of dyadic partitions of unity we refer to [BCD11, Proposition 2.10]. The Littlewood-Paley blocks are now defined as
$$\Delta_{-1}u=\mathcal{F}^{-1}(\chi\mathcal{F}u)\quad \Delta_{j}u=\mathcal{F}^{-1}(\theta(2^{-j}\cdot)\mathcal{F}u).$$
We point out that everything above and everything that follows can be applied to distributions on the torus (see [SW71]). More precisely,  Besov spaces on the torus with general indices $p,q\in[1,\infty]$ are defined as the completion of $C^\infty(\mathbb{T}^d)$ with respect to the norm
$$\|u\|_{B^\alpha_{p,q}}:=(\sum_{j\geq-1}(2^{j\alpha}\|\Delta_ju\|_{L^p(\mathbb{T}^d)})^q)^{1/q}.$$
 We  will need the following Besov embedding theorem on the torus (c.f. [GIP15, Lemma 41]):
\vskip.10in
 \th{Lemma A.1} i) Let $1\leq p_1\leq p_2\leq\infty$ and $1\leq q_1\leq q_2\leq\infty$, and let $\alpha\in\mathbb{R}$. Then $B^\alpha_{p_1,q_1}(\mathbb{T}^d)$ is continuously embedded in $B^{\alpha-d(1/p_1-1/p_2)}_{p_2,q_2}(\mathbb{T}^d)$.

 ii) (Besov embedding [Tri06, Chapter 6]) Let $\alpha_1<\alpha_2$, $1\leq p_1\leq p_2\leq \infty$,
and $1\leq q_1\leq q_2\leq \infty$. Then
$${{B}}^{\alpha_2}_{p_1,q_2}(\mathbb{T}^d)\subset {{B}}^{\alpha_1}_{p_1,q_1}(\mathbb{T}^d);\quad {{B}}^{\alpha_1}_{p_1,q_1}(\mathbb{T}^d)\subset {{B}}^{\alpha_1}_{p_1,q_2}(\mathbb{T}^d),\quad {{B}}^{\alpha_1}_{p_2,q_1}(\mathbb{T}^d)\subset {{B}}^{\alpha_1}_{p_1,q_1}(\mathbb{T}^d).$$

iii) ([MW15, Remarks 3.5, 3.6]) For $p>1$
$${{B}}^{0}_{p,1}(\mathbb{T}^d)\subset L^p\subset {{B}}^{0}_{p,\infty}(\mathbb{T}^d).$$

\vskip.10in

 Now we recall the following paraproduct introduced by Bony (see [Bon81]). In general, the product $fg$ of two distributions $f\in \mathcal{C}^\alpha, g\in \mathcal{C}^\beta$ is well defined if and only if $\alpha+\beta>0$. In terms of Littlewood-Paley blocks, the product $fg$ can be formally decomposed as
 $$fg=\sum_{j\geq-1}\sum_{i\geq-1}\Delta_if\Delta_jg=\pi_<(f,g)+\pi_0(f,g)+\pi_>(f,g),$$
 with $$\pi_<(f,g)=\pi_>(g,f)=\sum_{j\geq-1}\sum_{i<j-1}\Delta_if\Delta_jg, \quad\pi_0(f,g)=\sum_{|i-j|\leq1}\Delta_if\Delta_jg.$$

\vskip.10in
 The basic result about these bilinear operations is given by the following estimates:
\vskip.10in
 \th{Lemma A.2}(Paraproduct estimates, [Bon 81, MW16, Proposition A.7]) Let $\alpha,\beta\in \mathbb{R}$ and $p,p_1,p_2,q\in[1,\infty]$ be such that
$$ \frac{1}{p}=\frac{1}{p_1}+\frac{1}{p_2}.$$
 Then we have $$\|\pi_<(f,g)\|_{B^{\beta}_{p,q}}\lesssim \|f\|_{L^{p_1}}\|g\|_{B^{\beta}_{p_2,q}}\quad f\in L^{p_1}, g\in B^{\beta}_{p_2,q},$$
 and for $\alpha<0$, furthermore,
 $$\|\pi_<(f,g)\|_{B^{\alpha+\beta}_{p,q}}\lesssim \|f\|_{B^{\alpha}_{p_1,q}}\|g\|_{B^{\beta}_{p_2,q}}\quad f\in B^{\alpha}_{p_1,q}, g\in B^{\beta}_{p_2,q}.$$
 For $\alpha+\beta>0$ we have
 $$\|\pi_0(f,g)\|_{B^{\alpha+\beta}_{p,q}}\lesssim \|f\|_{B^{\alpha}_{p_1,q}}\|g\|_{B^{\beta}_{p_2,q}}\quad f\in B^{\alpha}_{p_1,q}, g\in B^{\beta}_{p_2,q}.$$

\vskip.10in
 The following basic commutator lemma is important for our use:
\vskip.10in
 \th{Lemma A.3}([GIP15, Lemma 5], [MW16, Proposition A.9]) Assume that $\alpha\in (0,1), \beta,\gamma\in \mathbb{R}$ and $p,p_1,p_2,p_3\in[1,\infty]$ are such that $$\alpha+\beta+\gamma>0,\quad \beta+\gamma<0,\quad \frac{1}{p}=\frac{1}{p_1}+\frac{1}{p_2}+\frac{1}{p_3}.$$ Then for smooth $f,g,h,$ the trilinear operator
 $$C(f,g,h)=\pi_0(\pi_<(f,g),h)-f\pi_0(g,h)$$ satisfies the bound
 $$\|C(f,g,h)\|_{B^{\alpha+\beta+\gamma}_{p,\infty}}\lesssim\|f\|_{B^\alpha_{p_1,\infty}}\|g\|_{B^\beta_{p_2,\infty}}\|h\|_{B^\gamma_{p_3,\infty}}.$$
 Thus, $C$ can be uniquely extended to a bounded trilinear operator from $B^\alpha_{p_1,\infty}\times B^\beta_{p_2,\infty} \times B^\gamma_{p_3,\infty}$ to $ B^{\alpha+\beta+\gamma}_{p,\infty}$.

\vskip.10in

Now we recall the following estimate for the heat semigroup $P_t:=e^{t\Delta}$.
\vskip.10in
\th{Lemma A.4}([GIP15, Lemma 47],[MW16, Proposition A.13] ) Let $u\in B^\alpha_{p,q}$ for some $\alpha\in \mathbb{R},p,q\in[1,\infty]$. Then for every $\delta\geq0$
$$\|P_tu\|_{B^{\alpha+\delta}_{p,q}}\lesssim t^{-\delta/2}\|u\|_{B^\alpha_{p,q}}.$$

\vskip.10in

\th{Lemma A.5} ([CC13, Lemma A.1]) Let $\alpha<1$ and $\beta\in \mathbb{R}$. Let $\varphi\in \mathcal{S}(\mathbb{R}^d)$, let $u\in \mathcal{C}^\alpha$, and $v\in \mathcal{C}^\beta$. Then for every $\varepsilon>0$
and every $\delta\geq-1$ we have
$$\|\varphi(\varepsilon D)\pi_<(u,v)-\pi_<(u,\varphi(\varepsilon D) v)\|_{\alpha+\beta+\delta}\lesssim \varepsilon^{-\delta}\|u\|_\alpha\|v\|_\beta.$$
where
$\varphi(D)u = \mathcal{F}^{-1}(\varphi \mathcal{F}u)$.
\vskip.10in

\th{Lemma A.6} ([CC13, Lemma 2.5], [MW16, Proposition A.13])  Let $u\in B_{p,q}^{\alpha+\delta}$ for some $\alpha\in \mathbb{R},0\leq\delta\leq2,p,q\in[1,\infty]$. Then for every $ t\geq0$
$$\|(P_t-I)u\|_{B^{\alpha}_{p,q}}\lesssim  t^{\delta/2}\|u\|_{B^{\alpha+\delta}_{p,q}}.$$
\vskip.20in
\no\textbf{Appendix B: Convergence of the stochastic terms}
\vskip.20in
We first recall the definition of the stochastic terms from [CC13] we use in the paper:
$$\aligned\includegraphics[height=0.5cm]{diri06.eps}:=&\lim_{\varepsilon\rightarrow0} \includegraphics[height=0.5cm]{diri07.eps} :=\lim_{\varepsilon\rightarrow0}(\includegraphics[height=0.5cm]{diri02.eps}^{2}-\bar{C}^{\varepsilon}_0),\\ \includegraphics[height=0.5cm]{diri11.eps} :=&\includegraphics[height=0.5cm]{diri02.eps}^3-3\bar{C}^{\varepsilon}_0\includegraphics[height=0.5cm]{diri02.eps},\\
\includegraphics[height=0.7cm]{diri14.eps}:=&\lim_{\varepsilon\rightarrow0}\includegraphics[height=0.5cm]{diri02.eps} \includegraphics[height=0.7cm]{diri04.eps}, \\ \includegraphics[height=0.7cm]{diri15.eps}:=&\lim_{\varepsilon\rightarrow0}\includegraphics[height=0.5cm]{diri07.eps} \diamond\includegraphics[height=0.7cm]{diri04.eps}:=\lim_{\varepsilon\rightarrow0}(\includegraphics[height=0.5cm]{diri07.eps} \includegraphics[height=0.7cm]{diri04.eps}-3(\bar{C}_{1}^\varepsilon+\bar{\varphi}^\varepsilon) \includegraphics[height=0.5cm]{diri02.eps}), \endaligned$$
$$\aligned
\includegraphics[height=0.5cm]{diri01.eps}\diamond(\includegraphics[height=0.7cm]{diri03.eps})^2
:=&\lim_{\varepsilon\rightarrow0}\includegraphics[height=0.5cm]{diri02.eps} (\includegraphics[height=0.7cm]{diri04.eps})^2
  \\\pi_{0,\diamond}(\includegraphics[height=0.7cm]{diri03.eps},\includegraphics[height=0.5cm]{diri06.eps} ) :=&\lim_{\varepsilon\rightarrow0}\pi_{0,\diamond}(\includegraphics[height=0.7cm]{diri04.eps},\includegraphics[height=0.5cm]{diri07.eps} ):=\lim_{\varepsilon\rightarrow0}(\pi_0(\includegraphics[height=0.7cm]{diri04.eps},\includegraphics[height=0.5cm]{diri07.eps})
 -3(\bar{C}_{1}^\varepsilon+\bar{\varphi}^\varepsilon) \includegraphics[height=0.5cm]{diri02.eps}),\\ \pi_{0,\diamond}(\includegraphics[height=0.7cm]{diri03.eps},\includegraphics[height=0.5cm]{diri01.eps} ) :=&\lim_{\varepsilon\rightarrow0}\pi_{0}(\includegraphics[height=0.7cm]{diri04.eps},\includegraphics[height=0.5cm]{diri02.eps} ),\\\pi_{0,\diamond}(\includegraphics[height=0.7cm]{diri08.eps},\includegraphics[height=0.5cm]{diri06.eps} ) :=&\lim_{\varepsilon\rightarrow0}\pi_{0,\diamond}(\includegraphics[height=0.7cm]{diri09.eps},\includegraphics[height=0.5cm]{diri07.eps} ):=\lim_{\varepsilon\rightarrow0}(\pi_0(\includegraphics[height=0.7cm]{diri09.eps},\includegraphics[height=0.5cm]{diri07.eps})
 -3(\bar{C}_{1}^\varepsilon+\bar{\varphi}^\varepsilon)),
 \\ \includegraphics[height=0.5cm]{diri07.eps}\diamond \bar{\Phi}_3^{\varepsilon}:=&\bar{\Phi}_3^{\varepsilon}\includegraphics[height=0.5cm]{diri07.eps}+3(\bar{C}_{1}^\varepsilon+\bar{\varphi}^\varepsilon) (-\includegraphics[height=0.7cm]{diri04.eps}+\bar{\Phi}_3^{\varepsilon}) .\endaligned$$
Here $\bar{C}_{0}^\varepsilon, \bar{C}_{1}^\varepsilon, \bar{\varphi}^\varepsilon$ are terms for renormalization and are defined in [CC13]. Here we do not recall the explicit formula of them since this is not used in our paper. The convergence above is in the corresponding  space (see (2.2)).
The convergence of $\delta \bar{C}_W^{\varepsilon}\rightarrow0$ can be obtained partially from [CC13] and a similar argument as in [CC13]. In this part we consider the convergence of $\int_0^T\includegraphics[height=0.5cm]{diri11.eps}ds$. We follow the notations from [GP17, Section 9].
 We represent the white noise in terms of its spatial Fourier transform. More precisely, let $E_0=\mathbb{Z}^3\backslash\{0\}$ and let $W(s,k)=\langle W(s),e_k\rangle$  and we view $W(s,k)$ as a Gaussian process on $\mathbb{R}\times E$ with covariance given by
 $$E\bigg[\int_{\mathbb{R}\times E_0}f(\eta) ˜W (d\eta) \int_{\mathbb{R}\times E_0}g(\eta') ˜W (d\eta')\bigg] =  \int_{R\times E_0}
g(\eta_1)f(\eta_{-1})d\eta_1,$$
where $\eta_a=(s_a,k_a)$, $s_{-a}=s_a, k_{-a}=-k_a$ and the measure $d\eta_a=ds_adk_a$ is the product of the Lebesgue measure $ds_a$ on $\mathbb{R}$ and of the counting measure $dk_a$ on $E_0$. Denote by $$\int_{(\mathbb{R}\times E_0)^n}f(\eta_{1...n})W(d\eta_{1...n})$$
  a generic element of the $n$-th chaos of $W$ on $\mathbb{R}\times E_0$. Recall that
 $$\int_0^t\includegraphics[height=0.5cm]{diri11.eps}d\sigma=2^{-3}\int_{(\mathbb{R}\times E)^3}e_{k_{[123]}}\int_0^tP_{\sigma-s_1}^\varepsilon(k_1)P_{\sigma-s_2}^\varepsilon(k_2)P_{\sigma-s_3}^\varepsilon(k_3)d\sigma W(d\eta_{123}).$$ Here $P_t^\varepsilon(k)=e^{-|k|^2t\pi^2}1_{\{t\geq0\}}g(\varepsilon k)$ and $k_{[123]}=k_1+k_2+k_3$.
By  a straightforward calculation we obtain that
$$\aligned &E|\Delta_q(\int_s^t(\bar{\Phi}_1^{\varepsilon_1})^{\diamond,3}d\sigma-\int_s^t(\bar{\Phi}_1^{\varepsilon_2})^{\diamond,3}d\sigma)|^2
\\\lesssim &\int_{(\mathbb{R}\times E)^2}\theta(2^{-q}k_{[123]})^2\bigg|\int_s^t[\Pi_{i=1}^3P_{\sigma-s_i}^{\varepsilon_1}(k_i)-\Pi_{i=1}^3P_{\sigma-s_i}^{\varepsilon_2}(k_i)]
d\sigma\bigg|^2d\eta_{123}\\\lesssim &(\varepsilon_1^\kappa+\varepsilon_2^\kappa)\int\theta(2^{-q}k_{[123]})^2\int_s^t\int_s^t\frac{e^{-\pi^2(|k_1|^2+|k_2|^2+|k_3|^2)|\sigma-\bar{\sigma}|}\sum_{i=1}^3|k_i|^\kappa}
{|k_1|^{2}|k_2|^{2}|k_3|^{2}}d\sigma d\bar{\sigma}dk_{123} \endaligned$$
$$\aligned
\lesssim &(\varepsilon_1^\kappa+\varepsilon_2^\kappa)\int\theta(2^{-q}k_{[123]})\frac{|t-s|\sum_{i=1}^3|k_i|^\kappa}
{|k_1|^{2}|k_2|^{2}|k_3|^{2}[|k_1|^2+|k_2|^2+|k_3|^2]}dk_{123}\\\lesssim &(\varepsilon_1^\kappa+\varepsilon_2^\kappa)\int_E\theta(2^{-q}k)\frac{|t-s|}{|k|^{2-\kappa}}dk
\lesssim(\varepsilon_1^\kappa+\varepsilon_2^\kappa)2^{q(1+\kappa)}|t-s|.\endaligned$$
Then by Gaussian hypercontractivity and  Lemma A.1  we obtain that  for any $\delta>0, p>1$, $\int_0^t\includegraphics[height=0.5cm]{diri11.eps}ds$ converges in $L^p(\Omega;C_T\mathcal{C}^{-\frac{1+\delta}{2}})$.

\vskip.20in
\no\textbf{Appendix C: Paracontrolled analysis for the solution to the lattice approximation}
\vskip.10in
In this appendix we recall paracontrolled analysis for the solution to (3.4) in [ZZ15].  To avoid confusion we do not use the graph notation for the lattice approximation in this paper. For the graph notation for $u^\varepsilon$  we refer to [ZZ15].  We define $$K^{\varepsilon}(t):=\int_0^tP_{t-s}^\varepsilon (u_1^{\varepsilon})^{\diamond,2}ds,\quad \tilde{K}^{\varepsilon}(t):=\int_0^t\tilde{P}_{t-s}^\varepsilon (u_1^{\varepsilon})^{\diamond,2}ds,$$
 and $$K_1^{\varepsilon}(t):=\int_0^tP_{t-s}^\varepsilon [e^{i_1i_2i_3}_N(u_1^{\varepsilon})^{\diamond,2}]ds,\quad \tilde{K}_1^{\varepsilon}(t):=\int_0^t\tilde{P}_{t-s}^\varepsilon [e^{i_1i_2i_3}_N(u_1^{\varepsilon})^{\diamond,2}]ds,$$
with $$\tilde{P}_t^\varepsilon:=\mathcal{F}^{-1}e^{-t|k|^2f(\varepsilon k)}\varphi_0(\varepsilon k)\mathcal{F},$$
where $\varphi_0$ is a smooth function and equals to $1$ on $\{|x|_\infty\leq 1\}$ with  $\textrm{supp}\varphi_0\subset\{|x|\leq 1.8\}$ and for $k=(k^1,k^2,k^3)\in\mathbb{R}^3$ $$f(k)=\frac{4}{|k|^2}(\sin^2\frac{k^1\pi}{2}+\sin^2\frac{k^2\pi}{2}+\sin^2\frac{k^3\pi}{2}).$$
Then we write the paracontrolled ansatz for the solution to (3.4) as follows: $$u_3^{\varepsilon}=-3P_N[\pi_<(u_2^{\varepsilon}+u_3^{\varepsilon},\tilde{K}^{\varepsilon}+\tilde{K}_1^{\varepsilon})]+u^{\varepsilon,\sharp}$$
 with $u^{\varepsilon,\sharp}(t)\in \mathcal{C}^{1+3\kappa}$.  Now we introduce the stochastic terms for the lattice approximation: for $T>0$
 $$\aligned C^\varepsilon_W(T):=&\sup_{t\in[0,T]}\big[\|u_1^{\varepsilon}\|_{-\frac{1}{2}-2\kappa}+\|(u_1^{\varepsilon})^{\diamond,2}\|_{-1-2\kappa}
 +\|u_2^{\varepsilon}\|_{\frac{1}{2}-2\kappa}+\|\pi_{0}( u_2^{\varepsilon},u_1^{\varepsilon})\|_{-2\kappa}\\&+\|\pi_{0,\diamond}( u_2^{\varepsilon},(u_1^{\varepsilon})^{\diamond,2})\|_{-\frac{1}{2}-2\kappa}+\|\pi_{0,\diamond}
(K^{\varepsilon},(u_1^{\varepsilon})^{\diamond,2})\|_{-2\kappa}\big]+\|u_2^{\varepsilon}\|_{C^{\frac{1}{8}}_T\mathcal{C}^{\frac{1}{4}-2\kappa}},\endaligned$$

$$\aligned E^\varepsilon_W(T):=&\sup_{t\in[0,T]}\big[\|(u_1^{\varepsilon})^{\diamond,2}e_N^{i_1i_2i_3}\|_{-1-2\kappa}+\|\pi_{0}( u_2^\varepsilon,e^{i_1i_2i_3}_N u^\varepsilon_1)\|_{-2\kappa}+\|\pi_{0,\diamond}(u_2^{\varepsilon},e^{i_1i_2i_3}_N(u_1^{\varepsilon})^{\diamond,2})\|_{-\frac{1}{2}-2\kappa}
\\&+\|\pi_0(K^\varepsilon,e_N^{i_1i_2i_3}(u_1^{\varepsilon})^{\diamond,2})\|_{-2\kappa}+\|\pi_0(K_1^\varepsilon,(u_1^{\varepsilon})^{\diamond,2})\|_{-2\kappa}
+\|\pi_{0,\diamond}(K_1^\varepsilon,e_N^{i_1i_2i_3}(u_1^{\varepsilon})^{\diamond,2})\|_{-2\kappa}\big]
,\endaligned$$
and$$\aligned\delta C_W^{\varepsilon}(T):=&\sup_{t\in[0,T]}\big[\|u_1^{\varepsilon}-\includegraphics[height=0.5cm]{diri01.eps}\|_{-\frac{1}{2}-2\kappa}
+\|(u_1^{\varepsilon})^{\diamond,2}-\includegraphics[height=0.5cm]{diri06.eps}\|_{-1-2\kappa}
+\|u_2^{\varepsilon}+\includegraphics[height=0.7cm]{diri03.eps}\|_{\frac{1}{2}-2\kappa}\\&+\|\pi_{0}( u_2^{\varepsilon},u_1^{\varepsilon})+\pi_{0,\diamond}( \includegraphics[height=0.7cm]{diri03.eps},\includegraphics[height=0.5cm]{diri01.eps})\|_{-2\kappa}+\|\pi_{0,\diamond}( u_2^{\varepsilon},(u_1^{\varepsilon})^{\diamond,2})+\pi_{0,\diamond}( \includegraphics[height=0.7cm]{diri03.eps},\includegraphics[height=0.5cm]{diri06.eps})\|_{-\frac{1}{2}-2\kappa}\\&+\|\pi_{0,\diamond}
(K^{\varepsilon},(u_1^{\varepsilon})^{\diamond,2})-\pi_{0,\diamond}
(\includegraphics[height=0.7cm]{diri08.eps},\includegraphics[height=0.5cm]{diri06.eps})\|_{-2\kappa}\big]
+\|u_2^{\varepsilon}+\includegraphics[height=0.7cm]{diri03.eps}\|_{C^{\frac{1}{8}}_T\mathcal{C}^{\frac{1}{4}-2\kappa}}.\endaligned$$
Here the terms containing $\diamond$ are renormlized terms defined in [ZZ15, Section 4]. Moreover, we introduce  the following operators
$$A_N^1(g,h)(f):=-\pi_0((I-P_N)\pi_<(f,P_Ng),h),$$
and $$A^2_N(g,h)(f):=\pi_0(P_N\pi_<(f,(P_{3N}-P_N)g),h).$$
Then we define
$$\aligned A_N(T):=&\|(A^1_N+A^2_N)(\tilde{K}^\varepsilon+\tilde{K}_1^\varepsilon,
(u_1^{\varepsilon})^{\diamond,2}+e_N^{i_1i_2i_3}(u_1^{\varepsilon})^{\diamond,2})\|_{C_TL(\mathcal{C}^{1-3\kappa},\mathcal{C}^{-\frac{1}{2}-5\kappa})}
\endaligned$$
and $$\aligned D_N(T):=&\sup_{t\in[0,T]}(\|-\pi_0((I-P_N)\pi_<(u_2^\varepsilon,K^\varepsilon+K_1^\varepsilon),
(u_1^{\varepsilon})^{\diamond,2}+e_N^{i_1i_2i_3}(u_1^{\varepsilon})^{\diamond,2})\\&+
\pi_0(P_N\pi_<(u_2^\varepsilon,(P_{3N}-P_N)(\tilde{K}^\varepsilon+\tilde{K}_1^\varepsilon)),
(u_1^{\varepsilon})^{\diamond,2}+e_N^{i_1i_2i_3}(u_1^{\varepsilon})^{\diamond,2})\|_{-\kappa}).\endaligned$$

By the calculations in [ZZ15] we obtain the following result.
\vskip.10in
\th{Proposition C.1} There exists $\kappa_1, C>0$ such that $$E[\delta C_W^{\varepsilon}(T)+ A_N(T)+ E^\varepsilon_W(T)+D_N(T)]\leq C \varepsilon^{\kappa_1}.$$
\vskip.10in
Moreover, by a similar argument as in [MW15, Lemma A.6] we obtain the following estimate on the extension operator defined in (3.2):
\vskip.10in
\th{Lemma C.2} Let $f$ be a function on $\Lambda_\varepsilon$. Then we have
$$\|\textrm{Ext} f\|_{L^{2n}(\mathbb{T}^3)}\lesssim N^{\frac{3}{2n}}\|f\|_{L^{2n}(\Lambda_\varepsilon)},$$
where the implicit constant depends on $n$.

\proof By (3.2) we have
$$\textrm{Ext} f(x)=\sum_{z\in \Lambda_\varepsilon}\frac{\varepsilon^3}{8}f(z)\Pi_{j=1}^3\frac{\sin\frac{\pi}{2}(2N+1)(x^j-z^j)}{\sin\frac{\pi}{2}(x^j-z^j)}.$$
Then we have
$$|\textrm{Ext} f(x)|^{2n}\lesssim \sum_{z\in \Lambda_\varepsilon}\frac{\varepsilon^3}{8}|f(z)|^{2n}[\sum_{z\in \Lambda_\varepsilon}\frac{\varepsilon^3}{8}\Pi_{j=1}^3\big{|}\frac{\sin\frac{\pi}{2}(2N+1)(x^j-z^j)}{\sin\frac{\pi}{2}(x^j-z^j)}\big{|}^{\frac{2n}{2n-1}}]^{2n-1}.$$
By the proof of [MW15, Lemma A.6] we obtain that
$$[\sum_{z\in \Lambda_\varepsilon}\frac{\varepsilon^3}{8}\Pi_{j=1}^3\big{|}\frac{\sin\frac{\pi}{2}(2N+1)(x^j-z^j)}{\sin\frac{\pi}{2}(x^j-z^j)}
\big{|}^{\frac{2n}{2n-1}}]^{2n-1}\lesssim N^3,$$
where the implicit constant does not depend on $x$, which implies the result.$\hfill\Box$

\vskip.10in
\no\textbf{Appendix D Symmetric quasi regular Dirichlet forms and Markov Processes}
\vskip.10in

In this section we recall some general Dirichlet form results from [MR92]. Let $E$ be a Hausdorff topological space, $m$ a $\sigma$-finite measure on $E$,
and let $\mathcal{B}$ the smallest $\sigma$-algebra of subsets of $E$ with respect to which all
continuous functions on $E$ are measurable.
Let $\mathcal{E}$ be a symmetric Dirichlet form acting in the real $L^2(m)$-space, i.e.
$\mathcal{E}$ is a positive, symmetric, bilinear, closed form with domain $D(\mathcal{E})$ dense
in $L^2(m)$, and such that $\mathcal{E}(\Phi(u),\Phi(u))\leq \mathcal{E}(u, u)$, for any $u\in D(\mathcal{E})$, where
$\Phi(t) = (0\vee t)\wedge1, t\in\mathbb{R}$. The latter condition is known to be equivalent
with the condition that the associated $C_0$-contraction semigroup $T_t, t\geq 0$,
is submarkovian (i.e. $0 \leq u\leq 1$ m-a.e. implies $0\leq T_tu \leq 1$ m-a.e., for
all $u\in L^2(m)$); association means that $\lim_{t\downarrow0}\frac{1}{t}\langle u- T_tu, v\rangle_{L^2(m)} = \mathcal{E}(u, v), \forall u, v\in D(\mathcal{E})$.
\vskip.10in

\th{Definition D.1} (cf. [MR92, Chap. IV, Defi. 3.1])
A symmetric Dirichlet form is called quasi-regular if the following holds:

i) There exists a sequence $(F_k)_{k\in\mathbb{N}}$ of compact subsets of $E$ such that $\cup_kD(\mathcal{E})_{F_k}$ is $\mathcal{E}_1^{1/2}$-dense in $D(\mathcal{E})$ (where $D(\mathcal{E})_{F_k}:= \{u \in D(\mathcal{E})|u =0 \textrm{ m-a.e. on } E-F_k\}$; $\mathcal{E}_1^{1/2}$ is the norm given by the scalar product in $L^2(m)$ defined by $\mathcal{E}_1$, where
$\mathcal{E}_1(u, v) := \mathcal{E}(u, v) + \langle u, v\rangle$,
$\langle , \rangle$ being the scalar product in $L^2(m)$. Such a sequence $(F_k)_{k\in\mathbb{N}}$
is called an $\mathcal{E}$-nest.

ii) There exists an $\mathcal{E}_1^{1/2}$-dense subset of $D(\mathcal{E})$ whose elements have $\mathcal{E}$-quasi continuous
$m$-versions. A real function $u$ on $E$ is called quasi continuous
when there exists an $\mathcal{E}$-nest $(F_k)$ s.t. $u$ restricted to $F_k$ is continuous.

iii) There exists $u_n\in D(\mathcal{E}), n\in\mathbb{N}$, with $\mathcal{E}$-quasi continuous $m$-versions $˜\tilde{u}_n$
and there exists an $\mathcal{E}$-exceptional subset $N$ of $E$ s.t. $\{˜\tilde{u}_n\}_{n\in\mathbb{N}}$ separates
the points of $E-N$. An $\mathcal{E}$-exceptional subset of $E$ is a subset
$N\subset \cap_k(E-F_k)$ for some $\mathcal{E}$-nest $(F_k)$.
\vskip.10in

To recall the main results in [MR92] we recall the definitions of a Markov process and a right process. Here we consider only Markov processes with life time $\infty$.
\vskip.10in

\th{Definition D.2} (cf. [MR92, Chap. IV Defi. 1.5]) A collection $\mathbf{M}:=(\Omega,\mathcal{M},(X_t)_{t\geq0},(P^z)_{z\in E})$ is called a Markov process (with state space $E$) if
it has the following properties.

i) There exists a filtration $(\mathcal{M}_t)$ on $(\Omega,\mathcal{M})$ such that $(X_t)_{t\geq0}$ is an $(\mathcal{M}_t)_{t\geq0}$ adapted stochastic process
 with state space $E$.

ii) For each $t\geq0$ there exists a shift operator $\theta_t:\Omega\rightarrow\Omega$ such that $X_s\circ\theta_t=X_{s+t}$ for all $s,t\geq0$

iii) $P^z, z\in E,$ are probability measures on $(\Omega,\mathcal{M})$ such that $z\mapsto P^z(\Gamma)$ is $\mathcal{B}(E)^*$-measurable for each $\Gamma\in\mathcal{M}$ resp. $\mathcal{B}(E)$-measurable if $\Gamma\in \sigma\{X_s|s\in[0,\infty)\}$, where $\mathcal{B}(E)^*:=\cap_{P\in\mathcal{P}(E)}\mathcal{B}^P(E)$ for $\mathcal{P}(E)$ denoting the family of all probability measures on $(E,\mathcal{B}(E))$ and $\mathcal{B}^P(E)$ denotes the completion of the $\sigma$-algebra $\mathcal{B}(E)$ w.r.t. a probability $P$.

iv) (Markov property) For all $A\in\mathcal{B}(E)$ and any $t,s\geq0$
$$P^z[X_{s+t}\in A|\mathcal{M}_s]=P^{X_s}[X_t\in A] \quad P^z-a.s., z\in E.$$

\vskip.10in

\th{Definition D.3} (cf. [MR92, Chap. IV Defi. 1.8]) Let $\mathbf{M}:=(\Omega,\mathcal{M},(X_t)_{t\geq0},(P^z)_{z\in E})$ be a Markov process with state space $E$ and corresponding filtration $(\mathcal{M}_t)$. $\mathbf{M}$ is called a right process if it has the following additional properties.

i) (Normal property) $P^z(X_0=z)=1$ for all $z\in E$.

ii) (Right continuity) For each $\omega\in\Omega$, $t\mapsto X_t(\omega)$ is right continuous on $[0,\infty)$.

iii) (Strong Markov property) $(\mathcal{M}_t)$ is right continuous and for every $(\mathcal{M}_t)$-stopping time $\sigma$ and every $\nu\in \mathcal{P}(E)$
$$P^\nu[X_{\sigma+t}\in A|\mathcal{M}_\sigma]=P^{X_\sigma}[X_t\in A]\quad P^\nu-a.s.$$
for all $A\in\mathcal{B}(E)$, $t\geq0$.

\vskip.10in

\th{Theorem D.4} (cf. [MR92, Chap. IV Thm 6.7]) Let $E$ be a metriable Lusin space. Then a Dirichlet form $(\mathcal{E},D(\mathcal{E}))$ on $L^2(E,m)$ is quasi-regular if and only if there exists a right process $\mathbf{M}$ associated with $(\mathcal{E},D(\mathcal{E}))$, i.e. the semigroup of $\mathbf{M}$ is an $m$-version of the semigroup associated with $(\mathcal{E},D(\mathcal{E}))$. In this case $\mathbf{M}$ is always properly associated with  $(\mathcal{E},D(\mathcal{E}))$.
\vskip.10in

\th{Remark D.5} The results in [MR92, Chap. IV] are more general and can be applied for general Hausdorff topological space and more general Markov process.
 Lusin spaces are enough for our use in this paper.


\begin{thebibliography}{99}
    \bibitem[ALZ06]{ } S. Albeverio, S. Liang and B. Zegarlinski, Remark on the integration by
parts formula for the $\Phi^4_3$-quantum field model. Infin. Dimens. Anal. Quantum Probab.
Relat. Top. 9, no. 1, (2006), 149-154.
\bibitem[AR91]{}S. Albeverio, M. R\"{o}ckner, Stochastic differential equations in infinite
dimensions: Solutions via Dirichlet forms, Probab. Theory Related Field 89 (1991) 347-386
\bibitem[BCD11]{} H. Bahouri, J.-Y. Chemin, R. Danchin,  Fourier analysis and nonlinear
partial differential equations, vol. 343 of Grundlehren der Mathematischen
Wissenschaften [Fundamental Principles of Mathematical Sciences]. Springer, Heidelberg,
2011.
\bibitem[BFS83]{} D. C. Brydges, J. Fr\"{o}hlich, and A. D. Sokal. A new proof of the existence
and nontriviality of the continuum $\Phi^4_2$ and $\Phi^4_3$ quantum field theories.
Comm. Math. Phys. 91, no. 2, (1983), 141-186.
\bibitem[BG97]{} L. Bertini, G. Giacomin,  Stochastic Burgers and KPZ equations from particle
systems. Comm. Math. Phys. 183, no. 3, (1997), 571-607.
\bibitem[Bon81]{} J.-M. Bony,  Calcul symbolique et propagation des singularit\'{e}s pour les \'{e}quations
aux d\'{e}riv\'{e}es partielles non lin\'{e}aires. Ann. Sci. \'{E}cole Norm. Sup. (4) 14, no. 2, (1981),
209-246.
\bibitem[Chung82]{}K. L. Chung, Lectures from Markov Processes to Brownian Motion, New York: Springer.
\bibitem[CC13]{}R\'{e}mi Catellier, Khalil Chouk, Paracontrolled Distributions and the 3-dimensional Stochastic Quantization Equation, arXiv:1310.6869
\bibitem[CFKZ08]{} Z.-Q. Chen, P. J. Fitzsimmons, K. Kuwae, and T.-S. Zhang, Stochastic calculus for symmetric Markov processes, The Annals of Probability
2008, Vol. 36, No. 3, 931-970
\bibitem[DD03]{}G. Da Prato, A. Debussche, Strong solutions to the stochastic quantization equations, Ann.
Probab., 31(4):1900-1916, (2003)
\bibitem[DZ02]{} G. Da Prato, J. Zabczyk, Second order partial differential equations in Hilbert spaces, Cambridge University Press (2002)
\bibitem[F95]{}P. J. Fitzsimmons, Even and odd continuous additive functionals. In Dirichlet Forms
and Stochastic Processes (Beijing, 1993) (Z.-M. Ma, M. R\"{o}ckner and J.-A. Yan, eds.)
139-154. de Gruyter, Berlin, (1995)
\bibitem[Fel74]{} J. Feldman, The $\Phi^4_3$ field theory in a finite volume. Comm. Math. Phys. 37, (1974),
93-120.
\bibitem[FOT94]{}M. Fukushima, Y. Oshima, and M. Takeda, Dirichlet Forms and Symmetric
Markov Processes. de Gruyter, Berlin (1994)
     \bibitem[GIP15]{} M. Gubinelli, P. Imkeller, N. Perkowski, Paracontrolled distributions and singular PDEs, Forum Math. Pi 3 no. 6(2015)
  \bibitem   [GJ87] {}J. Glimm, A. Jaffe, Quantum physics. Springer-Verlag, New York, second ed.,
1987. A functional integral point of view.
\bibitem[GJ13]{}P. Goncalves and M. Jara. Nonlinear Fluctuations of Weakly Asymmetric Interacting
Particle Systems. Archive for Rational Mechanics and Analysis, 212(2):597-644, 2013.
\bibitem[GJ13a]{} M. Gubinelli, M. Jara. “Regularization by noise and stochastic Burgers
equations.” Stochastic Partial Differential Equations: Analysis and Computations 1.2 (2013):
325-350.
\bibitem[GP15]{} M. Gubinelli, N. Perkowski, Energy solutions of KPZ are unique, http://arxiv.org/abs/1508.07764v1 (2015)
\bibitem[GP17]{}M. Gubinelli, N. Perkowski, KPZ reloaded, Communications in Mathematical Physics, 349(1):165-269, (2017)
\bibitem[GRS75]{}F. Guerra, J. Rosen, B. Simon: The $P(\Phi)_2$ Euclidean quantum field theory
as classical statistical mechanics. Ann. Math. 101,  11-259 (1975)
\bibitem[Hai13]{} M. Hairer, Solving the KPZ equation. Ann. of Math. (2) 178, no. 2, (2013), 559-664.
\bibitem[Hai14]{} M. Hairer, A theory of regularity structures. Invent. Math. 198(2), 269-504, (2014).
 \bibitem[HM15]{} M. Hairer,  K. Matetski, Discretisations of rough stochastic PDEs, http://arxiv.org/abs/1511.06937v1
 \bibitem[HM16]{}M. Hairer and J. Mattingly. The strong Feller property for singular stochastic PDEs. Preprint,
arXiv:1610.03415, 2016.
\bibitem[HP14]{} M. Hairer and \'{E}. Pardoux, A Wong-Zakai theorem for stochastic PDEs, arXiv:1409.3138
\bibitem[HS16]{} M. Hairer, H. Shen. The dynamical sine-Gordon model. Communications in Mathematical Physics, 2016, 341(3), 1-57
\bibitem[HS17]{}M. Hairer and P. Schoenbauer. The strong support theorem for singular stochastic PDEs. In
preparation
    \bibitem[JLM85]{} G. Jona-Lasinio and P. K. Mitter. On the stochastic quantization of field theory. Comm.
Math. Phys., 101(3):409-436, 1985.
    \bibitem[KPZ86]{} M. Kardar, G. Parisi, Y.-C. Zhang, Dynamic scaling of growing interfaces.
Phys. Rev. Lett. 56, no. 9, (1986), 889-892.
\bibitem[Kup16]{} A. Kupiainen. Renormalization group and stochastic PDE’s. Annales Henri Poincar\'{e}, 2016, 17(3):497-535.
\bibitem[Lyo98]{} T. J. Lyons, Differential equations driven by rough signals. Rev. Mat. Iberoamericana
14, no. 2, (1998), 215-310.
\bibitem[LR15]{}W. Liu, M. R\"{o}ckner. Stochastic Partial Differential Equations: An Introduction.
Springer, (2015).
\bibitem[MR92]{}Z. M. Ma and M. R\"{o}ckner, Introduction to the Theory of (Non-Symmetric) Dirichlet
Forms (Springer-Verlag, Berlin, Heidelberg, New York, 1992.
\bibitem[MW14]{} J.-C. Mourrat, H. Weber. Convergence of the two-dimensional dynamic Ising-Kac model
to $\Phi^4_2$. Preprint, arXiv:1410.1179, 2014.
\bibitem[MW15]{} J.-C. Mourrat, H. Weber, Global well-posedness of the dynamic $\Phi^4$ model in the plane, arXiv:1501.06191v1
\bibitem[MW16]{} J.-C. Mourrat, H. Weber, Global well-posedness of the dynamic $\Phi^4_3$ model on the torus, first edition, arXiv:1601.01234
\bibitem[MW17]{} J.-C. Mourrat, H. Weber, The dynamic $\Phi^4_3$ model comes down from infinity, arXiv:1601.01234
\bibitem[N85]{}S. Nakao, Stochastic calculus for continuous additive functionals of zero energy. Z.
Wahrsch. Verw. Gebiete 68 557-578, (1985)
\bibitem[P75]{}Y. M. Park, Lattice approximation of the $(\lambda\varphi^4-\mu\varphi)_3$ field theory in a finite volume, J. Math. Phys. 16, 1065 (1975);
    \bibitem[PR07]{} C. Prevot, M. R\"{o}ckner, A Concise Course on Stochastic Partial Differential Equations, Lecture Notes in Math., vol.1905, Springer, (2007)
 \bibitem[PW81]{} G. Parisi,  Y. S. Wu. Perturbation theory without gauge fixing. Sci. Sinica 24,
no. 4, (1981), 483–496.
\bibitem[RZZ15] { }M. R\"{o}ckner, R. Zhu, X. Zhu, Restricted Markov unqiueness for the
stochastic quantization of $P(\phi)_2$ and its
applications,  arXiv:1511.08030 (2015), to appear in Journal of functional analysis
\bibitem[RZZ16] { }M. R\"{o}ckner, R. Zhu, X. Zhu, Ergodicity for the stochastic quantization problems on the 2D-torus,  arXiv:1511.08030 (2015), to appear in Communication in Mathemathical physics,  arXiv:1606.02102, (2016)
\bibitem[SW71]{}E. M. Stein, G. L. Weiss, Introduction to Fourier Analysis on Euclidean Spaces,
Princeton University Press, 1971
\bibitem[Tri06]{}H. Triebel, Theory of function spaces III. Basel, Birkh\"{a}user, (2006)
\bibitem[W05]{}F.Y.Wang, Functional Inequalities, Markov Semigroup and Spectral Theory. Chinese Sciences
Press, Beijing (2005)
\bibitem[ZZ14]{} R. Zhu, X. Zhu, Approximating three-dimensional Navier-Stokes equations
driven by space-time white noise, arXiv preprint arXiv:1409.4864 (2014)
\bibitem[ZZ15] { }R. Zhu, X. Zhu, Lattice approximation to the dynamical $\Phi_3^4$ model,  arXiv:1508.05613, to appear in The annals of Probability
\bibitem[ZZ15a]{} R. Zhu, X. Zhu, Three-dimensional Navier-Stokes equations driven by space-time white noise, Journal of Differential Equations
, 259,  9, 5,  2015,  4443-4508
\end{thebibliography}
\end{document}